\documentstyle{amsppt}
\magnification=\magstep1
\topmatter
\title 
Quotients of K3 Surfaces Modulo Involutions
\endtitle
\author D. -Q. Zhang
\endauthor
\abstract
Let  $X$  be a K3 surface with an involution
$\sigma$  which has non-empty fixed locus  $X^{\sigma}$ 
and acts non-trivially on a non-zero holomorphic
$2$-form.  We shall construct
all such pairs  $(X, \sigma)$  in a canonical way, 
from some better known double coverings of 
log del Pezzo surfaces of index  $\le 2$  or rational elliptic 
surfaces, and construct the only family of each of the 
three extremal case where  $X^{\sigma}$  contains  $10$  
(maximum possible) curves. 
We also classify rational log Enriques surfaces of index  $2.$
\endabstract
\endtopmatter

\document
{\bf Introduction}

\par \vskip 1pc
Let  $X$  be a smooth projective K3 surface over the
complex number field  ${\bold C}.$  Let  ${\sigma}$
be an involution on  $X.$  This  $\sigma$
induces an action on a non-zero holomorphic 2-form
$\omega$  of  $X$  such that  $\sigma^*\omega = \pm \omega.$
If  $\sigma^* \omega = \omega$  then the quotient
space  $X/\sigma$  is again a K3 surface with at worst
type  $A_1$  Du Val singular points.

\par
We are interested in the case where  $\sigma^*\omega = -\omega.$
We want to determine the structure of  $S := X/\sigma$
and construct all such double
coverings  $\pi : X \rightarrow S,$  as canonical 
resolutions [H1, p. 48], of some better known double
coverings  $\pi_c : X_c \rightarrow S_c$ (Theorem 1).
$S_c$  is, except for the Enriques case
(Lemma 1.2(3)) and the log del Pezzo case (Lemma 3.2),
a rational elliptic surface and  $\pi_c$  ramifies
over one or two fibers.

\par
We refer to [R] for the comparison of advantages of 
two approaches : top-down and bottom-up.
For involution-actions on the K3 lattice,
we refer to [N1,2].

\par
In this paper we shall prove the following Theorems 1, 3 and 4 
and Corollary 5.  Theorem 1 below, in the cases of 
Lemmas 2.1 and 3.2, also tells the relations between 
rational log Enriques surfaces and log del Pezzo surfaces 
downstairs, which have been studied in 
[AN, Bl, D, MZ1,2, N3,4,5, Z1,2,3], 
and K3 surfaces with an involution upstairs.

\par \vskip 1.5pc
{\bf Theorem 1.} {\it Let  $X$  be a smooth projective
K$3$ surface with an involution  $\sigma$  such that
$\sigma^*\omega = -\omega$  for a non-zero holomorphic
$2$-form  $\omega.$  Let  $\pi : X \rightarrow S := X/\sigma$
be the quotient morphism.  Then the following two assertions
hold.} 

\par 
(1) {\it The fixed locus  $X^{\sigma}$
is empty if and only if  $S$  is an Enriques surface.}

\par
(2) {\it Suppose that  $X^{\sigma} \ne \phi.$  Then  $\pi$  
is the canonical resolution, in the sense of Horikawa {\rm [H1, p.48],} 
of a double covering  $\pi_c : X_c \rightarrow S_c$  
so that  $\pi$  and  $\pi_c$  are precisely
constructed in Lemma {\rm 2.4, 3.2, 4.1, 5.1} or {\rm 6.1.}  
In the cases of Lemmas {\rm 4.1} and {\rm 5.1}, we have  $\pi = \pi_c.$}

\par
{\it In particular, one has the following commutative 
diagram with  $f$  as the minimal resolution and  $g$  a
resolution}
$$X \overset{f} \to{\rightarrow} X_c$$
$$\pi \downarrow \hskip 2pc \downarrow \pi_c$$
$$S \overset{g} \to{\rightarrow} S_c.$$

\par
{\bf Remark 2.} (1) We call  $\pi_c : X_c \rightarrow S_c$  a
{\it canonical mapping model} of  $\pi : X \rightarrow S.$
This model is unique in the cases of Lemmas 3.2, 4.1, 5.1 and 6.1,
but not unique in the case of Lemma 2.4 (see Definition and Propositions
2.10, 3.5 and 6.4).

\par
(2) This theorem also shows the usefulness of Persson's list 
given in [P], because in all cases, except for the case of 
Lemma 3.2, the canonical mapping model is constructed 
from a rational relatively minimal elliptic fibration. 

\par \vskip 1.5pc
{\bf Theorem 3} (see Theorem 3' at the end of \S 1 for
the detailed version). {\it Let  $(X, \sigma)$  be the pair 
satisfying the hypothesis of Theorem {\rm 1}.  Let  
$X^{\sigma}$  be the fixed locus.
Then the following four assertions are true.}

\par
(1) {\it $X^{\sigma}$  is a disjoint union of  
$m$  smooth curves for some  $0 \le m \le 10$  and  $m$
can attain any value in this range.  If  $m = 10,$  
then  $(X, \sigma)$  has one of the following three types.}

\par
Type{\rm (Rat)}. {\it $X^{\sigma}$  is a union of  $10$
rational curves.}

\par
Type{\rm (Gn2)}. {\it $X^{\sigma}$  is a union of one
genus-two curve and  $9$  rational curves.}

\par
Type{\rm (Ell)}. {\it $X^{\sigma}$  is a union of one
elliptic curve and  $9$  rational curves.}

\par
(2) {\it There is, up to isomorphisms, only one pair
$(X, \sigma)$ of Type{\rm (Rat)}.}

\par
(3) {\it There is a family  ${\Cal X}_{ell} := 
\{(X_s, \sigma_s) | s \in {\bold P}^1\}$
of surface  $X_s$  and an involution  $\sigma_s$  on it
such that for  $s \ne \infty, 0, 1, s_0,$  the pair  
$(X_s, \sigma_s)$  is of Type{\rm (Ell)}.  Conversely, 
every pair of Type{\rm (Ell)} is isomorphic to  $(X_s, \sigma_s)$
for some  $s \ne \infty, 0, 1, s_0.$}

\par
(4) {\it There is a family  ${\Cal Y}_{gn2} = 
\{(Y_t, \sigma_t) | t \in {\bold P}^3\}$
of surface  $Y_t$  and an involution
$\sigma_t$  on it such that for a general  $t,$  the pair  
$(Y_t, \sigma_t)$  is of Type{\rm (Gn2)}.  
Conversely, every pair of Type{\rm (Gn2)} is isomorphic to  
$(Y_t, \sigma_t)$  for some  $t.$}

\par \vskip 1.5pc
{\bf Remark.} (1) Nikulin [N2] classified the configuration
of  $X^{\sigma}$  in the case  $\sigma^*|Pic X = id$.  
However, our argument and result are more geometrical,
even in this case.

\par
(2) We also prove in Theorem 3' that the Picard number  $\rho(X) \ge 18$  
in all three extremal cases.  We refer to [Mo] for K3 surfaces with
large Picard number.

\par
(3) Though we constructed a 3-dimensional family
${\Cal Y}_{gn2}$  of K3 surfaces of Type(Gn2) and 
Picard number $\ge 18$, we have not identified K3 
surfaces in the family which are isomorphic to each other,
and not all, I guess, of K3 surfaces of Picard
number  $\ge 18$  are included in this family.
Therefore, this dimension 3 may not give any
restriction on the dimensions of moduli
spaces of K3 surfaces with Picard number $\ge 18$.
 
\par \vskip 1.5pc
In view of Theorem 3' (3)(4), the extremal Type(Rat) is also,
with one or two smooth rational curves, contracted, 
the degeneration of other two extremal Types (Ell) and (Gn2).  
So {\it the most extremal pair} should be a right name for
the pair of Type(Rat) (see [OZ1,2]).
Also, Theorem 3'(4) supports the naming of extremal case in [AN]
for the surface  $S_{gn2}.$

\par
Our next Theorem 4 reduces, in the sense of the 3-column diagram
below, the classification of rational log Enriques surfaces of 
index  $2$  to that of pairs  $(X, \sigma)$  as in Theorem 1: 
the case of Lemma 2.4, 
and hence reduce further to that of double coverings of rational 
elliptic surfaces ramifying over one or two fibers (Lemma 2.4).

\par \vskip 1.5pc
{\bf Theorem 4.} {\it Let  $R$  be a rational log Enriques 
surface of index  $2$  with  
${\overline \pi} : W \rightarrow R$  as its canonical 
covering \rm{(Definition 1.7)}.}

\par
{\it Then there exists a $($smooth$)$  K$3$ surface  $X$  with an 
involution  ${\sigma}$  such that the fixed locus  $X^{\sigma}$  
is a disjoint union of  $n$ $(1 \le n \le 10)$ smooth rational 
curves, and that the quotient morphism  
$\pi : X \rightarrow S := X/\sigma$  is 
the canonical resolution of  ${\overline \pi}$  
\rm{(Definition 2.2)}.}

\par
{\it In particular, one has the following commutative diagram,
where  $\pi_c : X_c \rightarrow S_c$  is a canonical mapping
model of  $\pi$ {\rm(Definition 2.10)}, and where
$p, f$  are the minimal resolutions and  $q, g$
resolutions all precisely described in Lemmas \rm{2.1} and \rm{2.4}}
$$W \overset{p} \to{\leftarrow} X \overset{f} \to{\rightarrow} X_c$$
$${\overline \pi} \downarrow \hskip 0.9pc 
\pi \downarrow \hskip 1.2pc \downarrow \pi_c$$
$$R \overset{q} \to{\leftarrow} S \overset{g} \to{\rightarrow} S_c.$$

\par \vskip 1.5pc
{\bf Corollary 5.} {\it Let  $R$  be a rational log Enriques surface
of index  $2$  and let  $r$  be the number of singular points of
Cartier index  $2.$  Then the following two assertions are true.}

\par
(1) {\it One has  $1 \le r \le 10,$  and  $r$  can attain any 
value in this range.}

\par
(2) {\it There is, up to isomorphisms, only one rational
log Enriques surface of index  $2$  with  $r = 10$  
\rm{(see Example 2.8 for the construction of this surface)}.}

\par \vskip 1.5pc
{\bf Remark 6.}  Let  $R$  be the unique rational log
Enriques surface of Type  $A_{19}$  constructed
in [Z2, Example 3.2 and Theorem 3.6] (see [OZ1, Example 2] 
and Example 2.8 below for two more different constructions, 
of the same surface).  The uniqueness theorem [OZ1, Theorem 2]
is the answer to the question asked by Naruki, 
and Reid who also discussed it in [R, Example 6].
Let  ${\overline \pi} : W \rightarrow R$ 
be the canonical double covering (Definition 1.7)
and let  $\pi : X \rightarrow S = X/\sigma$ be the 
canonical resolution of  ${\overline \pi}$ (Definition 2.2).

\par
Then the pair  $(X, \sigma)$  here is isomorphic
to the pair of Type(Rat) in Theorem 3,
while the surface obtained from  $S$  by contracting
$\pi(X^{\sigma})$  into  $10$  cyclic-quotient singular
points of Brieskorn type  $C_{4,1},$  is isomorphic to 
the unique rational log Enriques surface of index  $2$
in Corollary 5(2) with  $r = 10$ (see proofs of 
Theorem 3 and Corollary 5 and Example 2.8). 
 
\par \vskip 1.5pc
The organisation of the paper is as follows.
The first section is on how the fixed locus  $X^{\sigma}$
sits in the surface  $X.$  Lemmas 1.5, 1.6 and 1.11,
which might be of general interest, describe how  $\sigma$  
acts on elliptic fibers and Dynkin diagrams.  

\par
$\S 2 \sim \S 6$  form the main
ingredients used in \S 7 to prove the theorems.
All three extremal pairs  $(X, \sigma)$  in Theorem 3 or 
Theorem 3' are precisely constructed in Example 2.8 and \S 7 
both by taking as  $S_c$ (see Theorem 1 for
the notation) the same smooth rational surface with a 
multiple-fiber free elliptic fibration  $\psi : S_c
\rightarrow {\bold P}^1$  having   $\{I_9, 3I_0\}$
as its set of singular fibers.  Such a pair
$(S_c, \psi)$  is unique up to isomorphisms (Lemma 7 in \S 7).
This  $S_c$  or the pair  $(S_c, \psi)$  will be called
{\it Persson's most extremal rational elliptic surface.}

\par 
As an application to [OZ1, Theorem 2], 
we show in Example 2.8 the following:

\par
Let  $S_{ci}$ ($i = 1,2$)  be a smooth rational 
surface with a multiple-fiber free 
elliptic fibration  $\psi_i : S_{ci} \rightarrow {\bold P}^1$  
which has two singular fibers  ${\overline F}_{i,j}$
($j = 1,2$) of Kodaira type II (setting  $n_{i,j} := 1$),
or III ($n_{i,j} := 2$), or IV ($n_{i,j} := 4$), or
$I_{n_{i,j}}$ ($n_{i,j} \ge 1$)  such that  
$n_{i,1} + n_{i,2} = 10.$  For instance, one can 
take Persson's pair as  $(S_{c_i}, \psi_i).$
Let  $g_i : S_i \rightarrow S_{ci}$  be a smooth 
blowing-up of points in  ${\overline F}_{i,j}$  so that 
$F_{i,j} := g_i^*{\overline F}_{i,j}$  fits respectively
Case ($\alpha$), or ($\beta$), or ($\gamma$), or ($\delta$)
in Lemma 1.5.  Then one has  $S_1 \cong S_2 \cong X/\sigma$
where  $(X, \sigma)$  is the unique pair of Type(Rat) in Theorem 3.

\par
The above result suggests that one can divide Persson's 
list in [P] of rational surfaces with an elliptic 
fibration into several classes so that any two members
in the same class can be transformed to each other by a
blowing-up succeeded by a blowing-down both similar to the
$g_i$  above. 

\par
In \S 7, we prove, as a corollary to Lemma 7,
that there is, up to isomorphisms, only one log del 
Pezzo surface with a type  
$A_8$  Du Val singular point as its only singular point.

\par \vskip 1.5pc
{\bf Acknowledgement.} The author would like to thank
the referee for very careful reading and suggestions
which make the paper clearer.

\par \vskip 2pc
{\bf \S 1. Preliminaries}

\par \vskip 1pc
The following Lemma 1.1 is also proved in [OZ1].

\par \vskip 1.5pc
{\bf Lemma 1.1.} {\it Let  $(X, \sigma)$  be 
the same as in Theorem $1$.  Then the following three 
assertions hold.}

\par
(1) {\it The fixed locus  $X^{\sigma}$  is 
a disjoint union of smooth curves.}

\par
(2) {\it Let $G_{1}$  and  $G_{2}$ be two $\sigma$-stable 
smooth rational curves with  $G_{1}.G_{2} = 1$. 
Then exactly one of  $G_i$  is  $\sigma$-fixed.}

\par
(3) {\it Let $G$ be a $\sigma$-stable but not $\sigma$-fixed 
smooth rational curve.  Then $\#(G \cap X^{\sigma}) = 2$.}

\par \vskip 1.5pc
{\it Proof.} (1)  By the hypothesis on  $\sigma,$ at any
$\sigma$-fixed point  $P,$  we have
$\sigma(x,y) = (x, -y)$  for suitable local coordinates at  $P.$
So  $P$  lies on the  $\sigma$-fixed curve  $y = 0,$
which is smooth.  

\par
(3) follows from Hurwitz's genus formula applied to
$\pi : G \rightarrow G/\pi.$

\par
(2)  Since  $G_1, \, G_2$  are  $\sigma$-stable, the 
intersection  $P := G_1 \cap G_2$  is  $\sigma$-fixed, i.e., 
$P \in X^{\sigma}.$  Note that a generator of the tangent
space to  $G_i$  at  $P,$  is an eigenvector of  $\sigma_*$
with respect to an eigenvalue  $\lambda_i.$
By the conditions on  $\sigma,$  one has  $\{\lambda_1, \lambda_2 \}$
$= \{1, -1\}$  as sets.  That is, for exactly one  $i \in \{1,2\},$
one has  $\lambda_i = 1$  and  $G_i$  is  $\sigma$-fixed.
This proves Lemma 1.1.

\par \vskip 1.5pc
{\bf Lemma 1.2.} {\it Let  $(X, \sigma)$  be as in Theorem $1$.
Let  $S := X/\sigma$  be the quotient space and  
$\pi : X \rightarrow S$  the quotient morphism.  
Then the following four assertions hold true.}

\par
(1) {\it $S$  is a smooth surface with irregularity  $q(S) = 0.$}

\par
(2) {\it Let  $X^{\sigma} = \coprod_{i=1}^m C_i,$  where  $C_i$
is a smooth irreducible curve.  Then  $\pi^*D_i = 2C_i$  
where  $D_i := \pi(C_i),$ $0 \sim \pi^*K_S + \sum_{i=1}^m C_i$  
and  $-2K_S \sim \sum_{i=1}^m D_i.$

\par
The quotient morphism  $\pi : X \rightarrow S$
coincides with the double covering 
\newline
$Spec \oplus_{i=0}^1 {\Cal O}(K_S)^{\otimes 2} \rightarrow S$
associated with the relation  
${\Cal O}(-K_S)^{\otimes 2} \cong {\Cal O}(\sum_{i=1}^m D_i).$}

\par
(3) {\it Suppose that  $X^{\sigma} = \phi.$  Then  $S$
is an Enriques surface and  $\pi : X \rightarrow S$
is the canonical unramified covering associated with the
relation  ${\Cal O}(K_S)^{\otimes 2} \cong {\Cal O}_S.$}

\par
(4) {\it Suppose that  $X^{\sigma} \ne \phi.$  Then  $S$
is a rational surface.}

\par \vskip 1.5pc
{\it Proof.}  (1) follows from Lemma 1.1 (1) and the fact that
$q(S) \le q(X) = 0.$

\par
(2) By the ramification formula  
$0 \sim K_X = \pi^*K_S + \sum_{i=1}^m C_i.$ 
This, together with the projection formula, implies
$-2K_S \sim \sum_{i=1}^m D_i.$  Now the second paragraph
of (2) follows by applying [H2, Lemma 2.1; H1, Lemma 4].

\par
(3) By the assumption on  $S$, we see that  $K_S$  is
not linearly equivalent to zero.  Now (3)
follows from (2) and (1).

\par
(4)  Since  $-2K_S \sim \sum_{i=1}^m D_i > 0,$
the pluri-genera  $P_n(S) = 0$  for all  $n \ge 1.$
Hence  $S$  is rational because  $q(S) = 0.$
This completes the proof of Lemma 1.2.

\par \vskip 1.5pc
{\bf Definition 1.3.}  Let  $T$  be a normal projective
surface with at worst quotient singular points.
For any  $t \in Sing T,$  there exists a ``small''
finite group  $G_t \subseteq GL_2({\bold C})$  such 
that  $(T, t)$  is isomorphic to  $({\bold C}^2/G_t, 0)$
analytically.  The index  $|G_t : G_t \cap SL_2({\bold C})|$
is called {\it the Cartier index} of  $t.$

\par
Note that  $t \in Sing T$  has Catier index  $1$  if and only
if  $t$  is a Du Val singular point.

\par
The Cartier index of a smooth point is defined as  $1.$
{\it The Cartier index}  $I = I(T)$  of  $T$  is defined
as the l.c.m. of Cartier indices of all points on  $T.$
Note that  $I$  is nothing but the smallest positive
integer such that the multiple  $IK_T$  of the 
canonical divisor  $K_T$  is a Cartier divisor.
 
\par \vskip 1.5pc
{\bf Lemma 1.4.} {\it Let  $T$  be a normal projective 
surface with at worst quotient singular points of index  $\le 2,$ 
and let  $t_i$  $(i = 1,2, \cdots, r; r \ge 0)$ 
be all singular points on  $T$  of index  $2.$
Let  $g_1 : S_1 \rightarrow T$  be the minimal resolution.
Then the following three assertions are true.}

\par
(1) {\it $t_i$  is cyclic of Brieskorn type  $C_{4n_i, 2n_i-1}$ 
$(n_i \ge 1)$.  Hence  $g_1^{-1}(t_i)$  is either a single  
$(-4)$-curve  $(n_i = 1)$,  or a linear chain of two  $(-3)$-curves as
tips and  $n_i-2$ $(n_i \ge 2)$ $(-2)$-curves as middle components:}
$$(-4), \,\,\, (-3)-(-2)-(-2)- \cdots -(-2)-(-2)-(-3).$$

\par
(2) {\it Let  $g_2 : S \rightarrow S_1$  be the smooth blowing-up
of all intersection points of  $g_1^{-1}(t_i)$  for
those  $i$  with  $n_i \ge 2.$  Let  $g = g_1 \circ g_2 : 
S \rightarrow T.$  Then  $g^{-1}(t_i)$  is either a single 
$(-4)$-curve  $D_{i,1}$ $(n_i = 1)$, or
a linear chain of  $n_i$  $(-4)$-curves  $D_{i,j}$ 
$(j = 1,2, \cdots, n_i ; n_i \ge 2)$ and  $n_i-1$ 
$(-1)$-curves  $H_{i,j}$ $(j = 1,2, \cdots, n_i-1)$:}
$$D_{i,1}, \,\,\, D_{i,1}-H_{i,1}-D_{i,2}-H_{i,2}- 
\cdots -D_{i,n_i-1}-H_{i,n_i-1}-D_{i,n_i}.$$

\par
(3) $K_S = g^*K_T - 1/2\sum_{i=1}^r \sum_{j=1}^{n_i} D_{i,j}.$

\par \vskip 1.5pc
{\it Proof.} Let  $B := g_1^{-1}(Sing T) = \sum_k B_k$  be 
the irreducible decomposition.  By the equivalence of quotient
singularity (resp. Du Val singularity) and log terminal singularity
(resp. canonical singularity) [Ka1, Cor. 1.9] (see also
[KMM] and [Ko]), there are rational numbers  $\alpha_i$  
with  $0 \le \alpha_i < 1$  such that
$$K_{S_1} = g_1^*K_T - \sum_k \alpha_k B_k.$$
Moreover, $\alpha_k = 0$  if and only  $g_1$
maps the connected component of  $B$  containing  $B_k$
into a Du Val singular point.  

\par
Since  $2K_T$  is Cartier by the hypothesis, $2\alpha_k$  is
an integer.  This, together with  $0 \le \alpha_k < 1,$
implies that  $\alpha_k = 0$  or  $1/2.$

\par
Now the second assertion in (1) is the consequence of
the calculation of the conditions 
$(K_{S_1} + \sum_k \alpha_k B_k) . B_t = 0$ (cf. the proof
of [Z2, Lemma 1.8]), whence the first assertion in (1) also follows.

\par
(2) follows from (1) and the construction of  $g_2,$  while 
(3) follows from the above ramification formula involving  $g_1.$
This proves Lemma 1.4.

\par \vskip 1.5pc
{\bf Lemma 1.5.} {\it Let  $X, \, \sigma$  and  
$\pi : X \rightarrow S := X/\sigma$  be as in Theorem $1$.
Suppose that there exists an elliptic fibration
$\varphi : X \rightarrow {\bold P}^1$  such that
$X^{\sigma}$  is contained in fibers.  Let  $E_i$  be a  
$\sigma$-stable singular fiber.  Then  
$Supp E_i = \sum_{j=1}^{n_i} C_{i,j} + \sum_j G_{i,j}$  
is one of the following five cases, where  $C_{i,j}$'s  are 
only  $\sigma$-fixed components contained in  $E_i.$}

\par
Case(A). {\it $E_i$  is of Kodaira type $IV$.  So  $E_i$  is a
union of three curves  $C_{i,1}, \, G_{i,1}, \, G_{i,2}$
which share the same point.  One has  $\sigma(G_{i,1}) = G_{i,2}.$}

\par
Case(B).  {\it $E_i$  is of Kodaira type  ${I_0}^*.$  Hence
$E_i$  is a tree with a central component  $G_{i,1}$  and four 
curves  $G_{i,2}, G_{i,3}, C_{i,1}, C_{i,2}$  sprouting out from  
$G_{i,1}.$  This  $G_{i,1}$  is  $\sigma$-stable
while  $\sigma(G_{i,2}) = G_{i,3}.$}
 
\par 
Case(C). {\it $E_i$  is of Kodaira type  $IV^*.$  So
$E_i$  is a tree with a central component  $C_{i,4}$  and three
twigs  $G_{i,j} + C_{i,j}$ $(j = 1,2,3)$ sprouting out from  $C_{i,4}$  
such that  $G_{i,j} . C_{i,4} = 1.$  Each  $G_{i,j}$  is  $\sigma$-stable.}

\par
Case(D). {\it $E_i$  is of Kodaira type  $I_{2n_i}.$  Hence
$E_i = C_{i,1}+G_{i,1}+C_{i,2}+G_{i,2}+ \cdots+ C_{i,n_i}+G_{i,n_i}$  
is a loop so that  $C_{i,j} . G_{i,j} = G_{i,j} . C_{i,j+1} = 1.$  
Here we set  $C_{i,n_i+1} = C_{i,1}, \,
G_{i,n_i+1} = G_{i,1}.$  Each  $G_{i,j}$  is  $\sigma$-stable.}

\par
Case(E).  {\it $E_i$  is of Kodaira type  $I_{2s_i}$ $(s_i \ge 1)$. 
So  $E_i = \sum_{j=1}^{2s_i} G_{i,j}$  is a loop so that
$G_{i,j} . G_{i,j+1} = 1.$  Here we set  $G_{i,2s_i+j} = G_{i,j}.$
One has  $\sigma(G_{i,j}) = G_{i,s_i+j}.$}
 
\par
{\it Finally, $\pi_*E_i$  fits one of the following cases
according to the type of  $E_i.$  Here  $\pi^*D_{i,j} = 2C_{i,j}$
and  $\sum_{j=1}^{n_i} D_{i,j}$  is a disjoint union of  $n_i$  
$(-4)$-curves  $D_{i,j}.$  One has also  $\pi^*H_{i,j} = G_{i,j},$  
except for Case$(\alpha)$ where  $\pi^*H_{i,1} = G_{i,1} + G_{i,2},$  
Case$(\beta)$ where  $\pi^*H_{i,2} = G_{i,2} + G_{i,3}$  
and Case$(\varepsilon)$ where $\pi^*H_{i,j} = G_{i,j} + G_{i,s_i+j}.$}

\par
Case$(\alpha)$. {\it $F_i := \pi_*E_i = 2H_{i,1} + D_{i,1}$  is a union of 
the touching  $(-1)$-curve  $H_{i,1}$  and  $(-4)$-curve  $D_{i,1}.$}

\par
Case$(\beta)$. {\it $F_i := \pi_*E_i = 4H_{i,1} + 2H_{i,2} + 
D_{i,1} + D_{i,2}$  
is a tree.  $H_{i,1}$  is a  $(-1)$-curve and also
the central component which meets the  $(-2)$-curve  
$H_{i,2}$  and the  $(-4)$-curves  $D_{i,1}, \, D_{i,2}.$}

\par
Case$(\gamma)$. {\it $F_i := \pi_*E_i = 3D_{i,4} + 4H_{i,1} + D_{i,1} + 
4H_{i,2} + D_{i,2} + 4H_{i,3} + D_{i,3}$  is a tree.  
$D_{i,4}$  is a  $(-4)$-curve and also the central component 
meeting three $(-1)$-curves  $H_{i,1}, \, H_{i,2}, \, H_{i,3},$  
while the  $(-4)$-curve  $D_{i,j}$ $(j = 1,2,3)$  is a tip 
component meeting only  $H_{i,j}.$}

\par
Case$(\delta)$. {\it $F_i := \pi_*E_i = D_{i,1} + 2H_{i,1} + D_{i,2} + 
2H_{i,2} + \cdots + D_{i,n_i} + 2H_{i,n_i}$  is a simple loop with  
$D_{i,j} . H_{i,j} = H_{i,j} . D_{i,j+1} = 1.$  Here  $D_{i,j}$
is a  $(-4)$-curve while  $H_{i,j}$  is a  $(-1)$-curve.}  

\par
Case$(\varepsilon)$. {\it $\pi_*E_i = 2F_i,$  where  
$F_i = \sum_{j=1}^{s_i} H_{i,j}$  is of Kodaira type  $I_{s_i}.$}
 
\par \vskip 1.5pc
{\it Proof.} Set  $E := E_i.$  By Lemma 1.1 and the hypothesis 
on  $\varphi,$  one obtains:

\par \vskip 1.5pc
{\bf Claim(1)}. (1.1)  If  $P$  is a point in  $E \cap X^{\sigma},$
then  $E \cap X^{\sigma}$  contains a smooth component 
through  $P.$

\par
(1.2) Suppose that  $G$  is a  $\sigma$-stable but not  
$\sigma$-fixed curve in  $E.$  Then  $E$  contains 
either one component  $C_1$  or two components  $C_1, \, C_2$  
of  $X^{\sigma}$  such that either  $\#(G \cap C_1) = 2$,  or  
$\#(G \cap C_i) = 1$  for both  $i = 1,2$, accordingly.

\par \vskip 1.5pc
By the classification theory, $E$  has one of the following
$10$  Kodaira types:

\par \vskip 1.5pc
(A) Type IV, \,\, (B) Type ${I_0}^*, \,\,$ (C) Type $IV^*, \,\,$
(DE) Type $I_m$ ($m \ge 2$),

\par
(F) $E$  is a nodal or cuspidal rational curve,

\par
(G) $E$  is a union of two touching smooth rational 
curves  $E_1, \, E_2,$  

\par
(H) Type ${I_n}^*$ ($n \ge 1$), and \,\, (I) Type $II^*$  or  $III^*.$

\par \vskip 1.5pc
If  $E$  is as in Case(F), then the singular point of  $E$
is  $\sigma$-fixed.  We reach a contradiction to Claim(1.1).

\par
If  $E$  is as in Case(G), then the common point of
$E_1$  and  $E_2$  is  $\sigma$-fixed.  By Claim(1.1) and Lemma 1 (1),
one may assume that  $E_1$  is  $\sigma$-fixed,
while  $E_2$  is  $\sigma$-stable but not  $\sigma$-fixed.  
This contradicts Claim(1.2) applied to  $G := E_2.$

\par
Let  $E$  be as in Case(H).  So  
$E = \sum_{i=1}^{n+1} E_i + E_{1,1}+E_{1,2} + E_{n+1,1}+E_{n+1,2}$  
consists of a linear chain  
$$E_1+E_2+ \cdots +E_n+E_{n+1}$$
and curves  $E_{i,1}, \, E_{i,2}$  sprouting out
from  $E_i$ ($i = 1$  and  $n+1$).
Note that either  $\sigma(E_i) = E_{n+2-i}$  for all  $i$  
or  $E_i$  is  $\sigma$-stable for all  $i.$

\par
In the first subcase, no  $E_i$  is  $\sigma$-fixed.
However, when  $n+1$  is odd (resp. even), the middle component  
$E_{(n+2)/2}$ (resp. the intersection  $E_{(n+1)/2} \cap E_{(n+1)/2+1}$)
is  $\sigma$-stable.  This contradicts Claim(1).

\par
In the second subcase, each  $E_i$ ($i = 1, n+1$)  is
not  $\sigma$-fixed by applying Lemma 1.1 (1) and Claim (1.2)
to  $G := E_{i,1}.$  Thus, at least one of  $E_{i,1}, E_{i,2},$ 
say  $E_{i,1},$ is  $\sigma$-fixed by Claim(1.2) applied 
to  $G := E_i.$  Then  $E_{i,2}$  is  $\sigma$-stable.
So  $E_i$  would have three  $\sigma$-fixed points
$E_i \cap E_{i,1}, \, E_i \cap E_{i,2}, \, E_i \cap E_j$
where  $j = i+1$ (resp. $i-1$)  if  $i = 1$ (resp. $i = n+1$).
Hence  $E_i$  must be  $\sigma$-fixed.  We reach a contradiction.

\par
Let  $E$  be as in Case(I).  Then  $E$  consists of a central
component  $R$  and three twigs  $T_i$  sprouting out from  $R.$
Note that  $R$  is  $\sigma$-stable.  But  $R$  is not  
$\sigma$-fixed by applying Lemma 1.1(1) and Claim(1.2) 
to  $G :=$ (the shortest twig among  $T_i$'s).  
Then, by Claim(1.2), $R$  meets a  $\sigma$-fixed curve in  
$T_i$  for  $i = 1$  and  $2$  say.  Hence all three twigs  $T_i$
are  $\sigma$-stable.  So  $R$  would have three  
$\sigma$-fixed points  $R \cap T_i,$  whence  $R$  must be 
$\sigma$-fixed.  This is a contradiction.

\par
By the same arguments as above, we can prove that Lemma 1.5 is true
if  $E$  is of type (A), (B), (C)  or  (DE).

\par \vskip 1.5pc
{\bf Lemma 1.6.} {\it Let  $(X, \sigma)$  be as in Theorem $1$.
Let  $\Gamma_i$  be a union of normal crossing smooth rational 
curves of Dynkin type  $A_m$ $(m \ge 1)$,  $D_m$ $(m \ge 4)$
or  $E_m$ $(m = 6,7,8)$.  Suppose that  $\Gamma_i$  is  $\sigma$-stable
and that every curve of  $X^{\sigma}$  is either containd in
$\Gamma_i$  or disjoint from  $\Gamma_i.$}

\par
{\it Then  $\Gamma_i$  is of type  $A_{2n_i-1}$  as follows:}
$$C_{i,1}-G_{i,1}-C_{i,2}-G_{i,2}- \cdots 
-C_{i,n_i-1}-G_{i,n_i-1}-C_{i,n_i}.$$
{\it Here  $C_{i,j}$  is  $\sigma$-fixed, while  $G_{i,j}$  is
$\sigma$-stable but not $\sigma$-fixed.}

\par \vskip 1.5pc
{\it Proof.}  The argument will be similar to Lemma 1.5.
In particular, the arguments for Lemma 1.5 Case(I) implies
that it is impossible that  $\Gamma_i$
is of type  $E_m$  or  $D_m.$  So  $\Gamma_i$  is of type  
$A_m$  as follows:
$$B_1-B_2- \cdots -B_m.$$
As in Lemma 1.5 Case(H), $\Gamma_i$  contains either a  $\sigma$-stable
curve or a  $\sigma$-fixed point.  Hence  $X^{\sigma} \cap \Gamma_i
\ne \phi$ (Lemma 1.1 (3)) and  $\Gamma_i$  contains a  $\sigma$-fixed
curve by the hypotheses of Lemma 1.6.  Thus  $\Gamma_i$  is 
component-wisely  $\sigma$-stable.  Now Lemma 1.6 follows from Lemma 1.1.

\par \vskip 1.5pc
{\bf Definition 1.7.} Let  $R$  be a normal projective
surface with at worst quotient singular points.

\par
(1) $R$  is {\it a log Enriques surface} if the irregularity
$q(R) = h^1(R, {\Cal O}_R) = 0$  and if
a positive multiple  $mK_R$  of the canonical
divisor  $K_R$  is linearly equivalent to zero.
{\it The index}  $I = I(R)$  is the smallest positive integer  
such that  $IK_R \sim 0,$ or equivalently the 
Cartier index in Definition 1.3 [Z2, Lemma 1.5].

\par
(2) The surface  $W := Spec\oplus_{i=0}^{I-1} {\Cal O}(-iK_R)$
or the natural quotient morphism   $\pi : W \rightarrow R,$
associated with the relation
${\Cal O}(K_R)^{\otimes I} \cong {\Cal O}_R,$
is called {\it the canonical covering} of  $R.$
This  $\pi$  is a Galois  ${\bold Z}/I{\bold Z}$-covering
such that  $W/({\bold Z}/I{\bold Z}) = R.$

\par
(3) In the sense of [OZ1], $R$  is called of {\it Type} 
$\alpha A_a + \delta D_d + \varepsilon E_e$
if the canonical covering  $W$  satisfies  
Sing $W = \alpha A_a + \delta D_d + \varepsilon E_e.$

\par \vskip 1.5pc
{\bf Remark 1.8.} (1) Note that  $W$  has at worst Du Val
singular points and  $K_W \sim 0.$  So  $W$  is either 
an abelian surface or a K3 surface with at worst Du Val
singular points.  Moreover,
$\pi : W \rightarrow R$  is unramified over  $R - Sing R.$

\par
(2) [Z2] classified the case where  $W$  is smooth.
The remaining cases of log Enriques surfaces
$R$  were dealt with in [Z3]; we 
proved there that there exists a crepant blowing up
$Q \rightarrow R$  with a new log Enriques surface  $Q$ 
of the same index such that the canonical cover  $V$  of
$Q$  has at worst type  $A_1$  Du Val singular points.
$Sing Q$  and  $\rho(Q)$  are tabulated there.

\par
(3) Blache [Bl] considers normal projective surfaces
with at worst log canonical singular points.
He also improves the upper bound of  $I$  to 
that  $I \le 21.$  Examples of log Enriques surfaces of 
all possible prime indices are given in [Z2].
  
\par
(4) Recently, we proved in [OZ1, 2] that
there is, upto isomorphisms, only one (resp. one or two) 
log Enriques surface(s) $R$  of Type  $A_{19},$ or  $D_{19},$ 
or  $D_{18}$  (resp. $A_{18}$).  
In the first case (resp. the last three cases),
the minimal resolution  $X$  of the canonical covering
$W$  of  $R,$  is the unique smooth K3 surface with 
discriminant of $Pic X$  equal to  $4$ (resp. $3$).

\par \vskip 1.5pc
{\bf Definition 1.9.}  Let  $T$  be a normal projective
surface with at worst quotient singular points.

\par
(1) $T$  is {\it a log del Pezzo surface} if the anti-canonical
divisor  $-K_T$  is a  ${\bold Q}$-ample divisor.

\par
(2)  A {\it Gorenstein} log del Pezzo surface is a log del Pezzo
surface of Cartier index  $1,$ equivalently, a log del Pezzo 
surface with at worst Du Val singular points.

\par \vskip 1.5pc
{\bf Remark 1.10.} (1) In view of [S], log del Pezzo 
surfaces are rational (see also [Z4, Lemma 1.1]
or [Z5, Lemma 1.3]).

\par
(2) Alexeev and Nikulin [AN] classified log del 
Pezzo surfaces of Cartier index  $2.$

\par
(3) [N3,4,5] gave upper bounds for the Picard number of
the minimal resolution of a log del Pezzo surface  $T$
in terms of the Cartier index or the l.c.m. of multiplicities of  $T.$ 

\par
(4) S. Keel and J. McKernan have announced in the 1995 Summer School of
Algebraic Geometry their affirmative answer to a conjecture of
Miyanishi, which says that the smooth part  $T - Sing T$ 
of a log del Pezzo surface  $T$  is rationally connected.

\par
(5) In [MZ 1, 2], Gorenstein log del Pezzo surfaces  
$T$  are classified by reducing to rank one or two cases, 
via a smooth blowing down.  In particular, it is proved 
there that the topological fundamental group  
$\pi_1(T - Sing T)$  is an abelian group of order  $\le 9.$

\par
(6) For an arbitrary log del Pezzo surface  $T,$  it is
proved in [GZ 1, 2] that the topological fundamental group  
$\pi_1(T - Sing T)$  is finite.  But this group may not
be abelian, in general;  see [Z1] for examples and also
the classification, when  $Y$  has Picard number one and
has at worst one rational triple and several 
Du Val singular points (97 types altogether).

\par \vskip 1.5pc
{\bf Lemma 1.11.} {\it Let  $(X,\sigma)$  be as in Lemma $1.1$
with  $X^{\sigma} \ne \phi.$
Suppose that there is an elliptic fibration
$\varphi : X \rightarrow {\bold P}^1$  such that  $X^{\sigma}$
is contained in fibers of  $\varphi$
and that there is a  $\sigma$-stable fiber
$E_1$  with  $E_1 \cap X^{\sigma} \ne \phi.$ 
Then  $S$  is a smooth rational surface
and the following three assertions are true.}

\par
(1) {\it $\sigma$  induces a permutation among fibers of  
$\varphi.$  There are exactly two  $\sigma$-stable fibers
$E_1, \, E_2$  of  $\varphi.$  We have  
$X^{\sigma} \subseteq Supp(E_1+E_2).$}

\par 
(2) {\it There is an elliptic fibration
$\psi : S \rightarrow {\bold P}^1$  with  $\pi_*E$
of a general fiber  $E$  of  $\varphi$  as a fiber of  $\psi.$ 
The pull back  $\pi^*\pi_*E$  is a disjoint union of
two smooth fibers  $E$  and  $\sigma(E)$  of  $\varphi.$}

\par
(3) {\it If  $mF$  is a fiber of  $\psi$  of multiplicity
$m \ge 2,$  then  $mF = \pi_*E_2.$  If  $F$  is not a minimal fiber 
of  $\psi,$ then  $F = \pi_*E_i$  for  $i = 1$  or  $2.$}

\par \vskip 1.5pc
{\it Proof.}  For every  $\sigma$-stable fiber  $E_i$  
(e.g. when  $i = 1$), Lemma 1.5 implies that if  $E_i$
contains (resp. does not contain) a component of  $X^{\sigma}$ 
then  $E_i$  fits one of Cases (A), (B), (C) and (D) 
(resp. is an elliptic curve or fits Case(E)) there
and  $\pi_*E_i$  is not a multiple fiber (resp.
$\pi_*E_i$  has multiplicity two).
In particular, (3) is a consequence of (1) and Lemma 1.5, 
together with the multiple-fiber freeness of an elliptic 
fibration on a K3 surface.

\par
Let  $E$  be a general fiber of  $\varphi.$
Then  $\sigma(E) \sim \sigma(E_1) = E_1.$  So  $\varphi$
induces a permutation among fibers of  $\varphi,$  and also
an automorphism on the base curve  ${\bold P}^1$  of  $\varphi.$

\par \vskip 1.5pc
{\bf Claim(1).}  $\sigma$  is not the identity automorphism
of  ${\bold P}^1.$  

\par \vskip 1.5pc
Supppse the contrary that Claim(1) is false.  Then 
every fiber of  $\varphi$  is  $\sigma$-stable.
Moreover, for a general fiber  $E$  of  $\varphi,$
the map  $\pi : E \rightarrow \pi(E)$  is
an etale covering of degree two.
For two general fibers  $E, \, E'$  of  $\varphi,$
one has  $2\pi(E) = \pi_*E \sim \pi_*E' = 2\pi(E').$
Hence  $\pi(E) \sim \pi(E')$  because  $S$  is a smooth rational
surface (Lemma 1.2).  Thus there is an elliptic fibration
$h : S \rightarrow {\bold P}^1$  with  $\pi(E)$  as
a general fiber.

\par
Now, $2\pi(E) = \pi_*E \sim \pi_*E_1.$
Since the g.c.d. of coefficients in  $\pi_*E_1$
is  $1$ (Lemma 1.5), one has
$m\pi_*E_1 \sim \pi(E)$  where  $m$ ($\ge 1$) is 
the multiplicity.  We reach a contradiction.  
This proves Claim(1).

\par \vskip 1.5pc
By Claim(1), $\sigma$  is an automorphism of order 2
on the base curve  ${\bold P}^1$  of  $\varphi.$
So, $\sigma$  has exactly two  $\sigma$-fixed points
by the Hurwitz genus formula for the covering
${\bold P}^1 \rightarrow {\bold P}^1/\sigma.$
Thus, there are exactly two
$\sigma$-stable fibers  $E_1, \, E_2$  of  $\varphi.$
This proves (1).

\par
Now we prove (2).  Let  $E$  be a general fiber of
$\varphi.$  By (1), the map  $\pi : E \rightarrow \pi(E)$
$= \pi_*E$  is an isomorphism.  For two general fibers  
$E, \, E'$  we have  $\pi_*E \sim \pi_*E'.$  So there is 
an elliptic fibration  $\psi : S \rightarrow {\bold P}^1$  with
$\pi_*E$  as a fiber and satisfying the conditions in (2).
This proves Lemma 1.11.

\par \vskip 1.5pc
In the subsequent sections, we shall also prove Theorem 3' below
which is stronger than Theorem 3 in the Introduction.

\par \vskip 1.5pc
{\bf Theorem 3'.} {\it Let  $(X, \sigma)$  be the pair 
satisfying the hypothesis of Theorem $1$ in the Introduction.  
Let  $X^{\sigma}$  be the fixed locus and  $\rho(X)$  the Picard 
number.  Then the following four assertions are true.}

\par
(1) {\it $X^{\sigma}$  is a disjoint union of  
$m$  smooth curves for some  $0 \le m \le 10$  and  $m$
can attain any value in this range.  If  $m = 10,$  
then  $(X, \sigma)$  has one of the following three types.}

\par
Type(Rat). {\it $X^{\sigma}$  is a union of  $10$
smooth rational curves.  One has  $\rho(X) = 20.$}

\par
Type(Gn2). {\it $X^{\sigma}$  is a union of one
genus-two curve and  $9$  smooth rational curves.
One has  $\rho(X) \ge 18.$}

\par
Type(Ell). {\it $X^{\sigma}$  is a union of one
elliptic curve  $E$  and  $9$  smooth rational curves.
One has  $\rho(X) \ge 19.$}

\par
(2) {\it There is, up to isomorphisms, only one pair
$(X, \sigma)$ of Type{\rm (Rat)}.  Such a pair is called
Shioda-Inose's pair in {\rm [OZ1, Example 2]
(see [Z2, Example 3.2]
and Example 2.8 below for two more constructions, 
of the same pair)}.  In particular, $X$  is isomorphic to 
the unique  K3 surface with the discriminant of  $Pic X$  
equal to  $4$ {\rm (see [V] for $Aut X$)}.}

\par
(3) {\it There is a family  ${\Cal X}_{ell} := 
\{(X_s, \sigma_s) | s \in {\bold P}^1\}$
of surface  $X_s$  and an involution  $\sigma_s$  on it,
satisfying the following four assertions
{\rm (see the proof for the construction)}.}

\par
(3i) {\it $X_s/\sigma_s$  is equal to  $S_{ell}$  
for all  $s$  with a smooth rational surface  $S_{ell},$ 
which is independent of  $s$  and given in the proof.
Let  $\pi : X_s \rightarrow S_{ell}$  be the quotient
morphism.}

\par
(3ii) {\it For  $s = \infty,\,\,$  
$X_{\infty} = S_{ell} \coprod S_{ell}$  is a disjoint 
union and  $\sigma_{\infty}$  is the order switching.}

\par
(3iii) {\it For three fixed distinct points  $s = 0, 1, s_0,$ 
the pair  $(Y_s, \sigma_s)$  is obtained from 
Shioda-Inose's pair  $(X, \sigma)$  so
that  $X \rightarrow X_s$  is the contraction of a 
$($ $\sigma$-stable$)$ smooth rational curve on  $X$  
meeting at two distinct points with a  $\sigma$-fixed curve
which is hence mapped to a  $1$-nodal curve  $E_s$  on  
$X_s$ $($the node $=$ {\rm Sing} $X_s$ $)$, and that  $\sigma_s$  is 
induced from  $\sigma.$  The fixed locus  $X_s^{\sigma_s}$
is a disjoint union of  $E_s$  and  $9$  smooth rational
curves.}

\par
(3iv) {\it For  $s \ne \infty, 0, 1, s_0,$  the pair  
$(X_s, \sigma_s)$  is of Type{\rm (Ell)}.  Conversely, 
every pair of Type(Ell) is isomorphic to  $(X_s, \sigma_s)$
for some  $s \ne \infty, 0, 1, s_0.$}

\par
(4) {\it There is a family  ${\Cal Y}_{gn2} = 
\{(Y_t, \sigma_t) | t \in {\bold P}^3\}$
of surface  $Y_t$  and an involution
$\sigma_t$  on it, satisfying the following three assertions 
{\rm (see the proof for the construction)}.}

\par
(4i) {\it $Y_t/\sigma_t$  is equal to  $S_{gn2}$  
for all  $t$  with a smooth rational surface  $S_{gn2},$  
which is independent of  $t$  and given in the proof.}

\par
(4ii) {\it For a general  $t,$  the pair  
$(Y_t, \sigma_t)$  is of Type{\rm (Gn2)}.  
Conversely, every pair of Type{\rm (Gn2)} is isomorphic to  
$(Y_t, \sigma_t)$  for some  $t.$}

\par
(4iii) {\it There is a straight projective line  $B$ 
in  ${\bold P}^3$  such that the subset  
${\overline {\Cal X}}_{ell} :=$ 
$\{({\overline X}_s, \sigma_s)$ 
$:= (Y_s, \sigma_s) | s \in B = {\bold P}^1 \}$  is obtained 
from the family  ${\Cal X}_{ell},$ in the following way:}

\par
{\it For  $s = \infty,$  ${\overline X}_s =
S_{gn2} \coprod S_{gn2}$  is a disjoint union
and  $\sigma_{\infty}$  is the order switching.
$S_{ell} \rightarrow S_{gn2}$  is the smooth blowing-down
of a  $(-1)$-curve  $H_0.$}
\par 
{\it For  $s \ne \infty,$  $X_s \rightarrow {\overline X}_s$
is the contraction of the  $\sigma_s$-stable
divisor  $\pi^*H_0$  on  $X_s$  meeting (with intersection
two) with  $X_s^{\sigma_s}$'s  only 
arithmetic-genus  $1$  curve  $E_s$  which is hence
mapped to an arithmetic genus  $2$  curve
${\overline E}_s$  with {\rm Sing} ${\overline E}_s =$ {\rm Sing}
${\overline X}_s$, and that the  $\sigma_s$  on  
${\overline X}_s$  is induced from the  $\sigma_s$  on  $X_s.$ 
The set  ${\overline X}_s^{\sigma_s}$
is a disjoint union of  ${\overline E}_s$  and  $9$  
smooth rational curves.

\par
To be precise, for some  $s_1 \in {\bold P}^1 - \{\infty, 0, 1, s_0\}$, 
$\pi^*H_0$  is a linear chain of two  $(-2)$-curves and  $E_{s_1}$
is a smooth elliptic curve through (transversally) the intersection
of  $\pi^*H_0$, and for all  $s \in {\bold P}^1 - \{\infty, s_1\}$,
$\pi^*H_0$  is a  $(-2)$-curve intersecting  $E_s$  at its
two distinct smooth points.}

\par \vskip 2pc
{\bf \S 2.  Rational log Enriques surfaces of index 2}

\par \vskip 1pc
Let  $R$  be a rational log Enriques surface of index  $2.$ 
We shall use the notation in Lemma 1.4 :
$$q_1 (= g_1) : S_1 \rightarrow R (= T), \,\,
q_2 (= g_2) : S \rightarrow S_1, \,\, q = q_1 \circ q_2, \,\,
K_S = q^*K_R - 1/2 \sum_{i=1}^r \sum_{j=1}^{n_i} D_{i,j}.$$

\par
Now the following relation
$$(2.1.1) \,\,\,\,\, {\Cal O}(-K_R)^{\otimes 2} \cong {\Cal O}_R$$
induces a relation:
$$(2.1.2) \,\,\,\,\, {\Cal O}(-K_S)^{\otimes 2} \cong 
{\Cal O}(\sum_{i=1}^r \sum_{j=1}^{n_i} D_{i,j}).$$

\par
Let  ${\overline \pi} : W := Spec \oplus_{i=0}^1$ 
${\Cal O}(iK_R) \rightarrow R$  and  
$\pi : X := Spec \oplus_{i=0}^1 {\Cal O}(iK_S) \rightarrow S$
be the double coverings associated with the relations 
(2.1.1) and (2.1.2), respectively.  Then both  $Gal(W/R)$  
and  $Gal(X/S)$  are isomorphic to the same 
cyclic group  $<\sigma>$  of order  $2$  such that  
$W/\sigma = R$  and  $X/\sigma = S.$

\par \vskip 1.5pc
{\bf Lemma 2.1.} {\it Let  $R$  be a rational log Enriques 
surface of index  $2.$  Let  ${\overline \pi} : W \rightarrow R =$
$W/\sigma$  and  $\pi : X \rightarrow S = X/\sigma$ 
be as above.  Then the following five assertions hold true.}

\par
(1) {\it One has  $\pi^*D_{i,j} = 2C_{i,j}$  for a smooth
rational curve  $C_{i,j}.$  The fixed locus
$X^{\sigma}$  is a disjoint union of  $n$ 
smooth rational curves  $C_{i,j},$  where  
$n := \sum_{i=1}^r n_i.$}

\par
(2) {\it $X$  is a smooth K$3$ surface.  The involution
$\sigma$  on  $X$  satisfies  $\sigma^*\omega = -\omega$  
for a non-zero holomorphic $2$-form  $\omega$  on  $X.$}

\par
(3) {\it There exists a  $\sigma$-equivariant birational
morphism  $p : X \rightarrow W$  which induces the following 
commutative diagram with  $p$  as the minimal resolution and  
$q$  a resolution}
$$X \overset{p} \to{\rightarrow} W$$
$$\pi \downarrow \hskip 2pc \downarrow {\overline \pi}$$
$$S \overset{q} \to {\rightarrow} R.$$

\par
(4) {\it There are  $n-1$  smooth rational curves  $G_j$  such that
$\Gamma := \sum_{i,j} C_{i,j} + \sum_j G_j$  is a component-wisely 
${\sigma}$-stable linear chain of length  $2n-1$  as follows: 
$$C_1-G_1-C_2-G_2- \cdots -C_{n-1}-G_{n-1}-C_n.$$
Here  $\{C_{i,j}\} = \{C_k\}$  as sets.
In particular, rank $Pic X \ge 2n$ and  $1 \le r \le n \le 10.$}

\par
{\it Let  $u_2 : X \rightarrow X_2$  be the contraction 
of  $\Gamma$  into a type  $A_{2n-1}$  Du Val singular point.  
Then  $S_2 := X_2/\sigma$  is a rational log Enriques 
surface of index  $2$  and Type  $A_{2n-1}$ {\rm (Def. 1.7)}.}

\par
(5) {\it Suppose in additoin that  $n = 10.$  
Then  $S_2$  is isomorphic to the unique rational 
log Enriques surface of Type  $A_{19}$
{\rm (see [OZ1, Theorem 2 and Example 2], [Z2, Example 3.2] and
Example 2.8 below for the proof of the uniqueness and 
three different constructions, of the same surface)}.
In particular, $X$  is isomorphic to the unique  K$3$ surface 
with the discriminant of  $Pic X$  equal to  $4.$}

\par \vskip 1.5pc
{\it Proof.} (1) Note that  $X$  is a smooth surface and 
$\pi$  is ramified exactly over the disjoint union
$\coprod D_{i,j}.$  So  $\pi^*D_{i,j} = 2C_{i,j}$  
for a smooth rational  curve  $C_{i,j}.$  Now (1) follows.

\par
(2) By the ramification formula, one has  $K_X \sim 0.$
Thus  $X$  is a smooth K3 surface because  $X$
is a double covering of a rational surface  (cf. [TY, Theorem 0.1]).

\par
Since  $\sigma^2 = id,$ one has  $\sigma^*\omega = \pm \omega.$
If  $\sigma^*\omega = \omega,$  then  $\sigma$  acts 
trivially on  $H^2(X, {\Cal O}_X) = {\bold C} {\overline \omega}.$
Hence  $H^2(S, {\Cal O}_S) = H^2(X, {\Cal O}_X)^{\sigma} \cong
{\bold C}.$  This contradicts the rationality of  $S$ 
because  $h^2(S, {\Cal O}_S) = h^0(S, K_S) = 0.$

\par
(3) follows from the constructions of  ${\overline \pi}, \, \pi.$

\par
(4) Let  $v_3 : S \rightarrow S_3$  be the 
contraction of  $n$  $(-4)$-curves  $D_{i,j}$  into  
$n$  cyclic quotient singular points of Brieskorn type  
$C_{4,1}.$  By the relation (2.1.2), one has
${\Cal O}(-K_{S_3})^{\otimes 2} \cong {\Cal O}_{S_3}.$
Hence, $S_3$  is a rational log Enriques
surface of index  $2.$  Now applying  [Z2, Theorem 3.6]  
to  ${\overline V} = S_3$  and  $(V, D) := (S, \sum_{i,j} D_{i,j}),$  
one sees that there are  $n-1$  $(-1)$-curves 
$H_j$  such that  $\Delta := \sum_{i,j} D_{i,j} + \sum_j H_j$
is a linear chain as follows, where  $\{D_{i,j}\} = \{D_k\}$  as sets
$$D_1-H_1-D_2-H_2- \cdots -D_{n-1}-H_{n-1}-D_n.$$
The pull back by  $\pi^*$  of this linear chain gives the 
linear chain  $\Gamma$  in (4), where  $G_j := \pi^*H_j.$  
Now the first paragraph of (4) is clear, while the second
paragraph follows from the observation that  $X_2$  
is a single singular point  $u_2(\Gamma).$

\par 
(5) follows from (4) and [OZ1, Theorem 2 and Example 2].

\par \vskip 1.5pc
{\bf Definition 2.2.}  In Lemma 2.1, 
the double covering  $\pi : X \rightarrow S$  is called 
{\it the canonical resolution} of the double covering  
${\overline \pi} : W \rightarrow R.$

\par \vskip 1.5pc
The following Lemma 2.3 is the converse to Lemma 2.1.

\par \vskip 1.5pc
{\bf Lemma 2.3} {\it Let  $(X, \sigma)$  be as in Theorem $1$.
Assume further that  $X^{\sigma}$  is a disjoint union
of  $n$  $(n \ge 1)$  smooth rational curves.

\par
Then the pair  $(X, \sigma)$  can be realized in the way of 
Lemma $2.1$, from a rational log Enriques surface  $R$  of 
index  $2$  and Type  $A_{2n-1}$  with a cyclic quotient 
singularity of Brieskorn type  $C_{4n, 2n-1}$  as its only 
singular point, i.e., the quotient morphism  
$\pi : X \rightarrow S := X/\sigma$  
is the canonical resolution of the double covering
${\overline \pi} : W \rightarrow R$  there.}

\par \vskip 1.5pc
{\it Proof.} By Lemma 1.2 (2), the quotient morphism 
$\pi : X \rightarrow S = X/\sigma$  coincides with
the double covering associated with the relation  (2.1.2)  
where  $\sum_{i,j} D_{i,j} = \pi(X^{\sigma}).$
By the proof of Lemma 2.1 (4), there are  $n-1$  $(-1)$-curves
$H_j$  on  $S$  such that  $\Delta :=$ 
$\sum_{i=1}^n D_i + \sum_{j=1}^{n-1} H_j$
is a linear chain as shown there, where  
$\sum_{i,j} D_{i,j} = \sum_{k=1}^n D_i.$

\par
Let  $q : S \rightarrow R$  be the contraction of  $\Delta$  
into a cyclic-quotient singular point of Brieskorn type  
$C_{4n, 2n-1}.$  Then our relation (2.1.2) induces the 
relation (2.1.1).  Hence  $R$  is a ratoinal log Enriques 
surface of index  $2$ (Lemma 1.2 (4)).

\par
Clearly,  $\pi$  is the canonical resolution of the
canonical double covering  ${\overline \pi} : W \rightarrow R$
associated with the relation (2.1.1), because
our contraction  $q : S \rightarrow R$  coincides with
the map  $q (= g) : S \rightarrow R (= T)$  constructed 
in Lemma 1.4.  This proves Lemma 2.3.

\par \vskip 1.5pc
We now construct the pair  $(X, \sigma)$  of Lemma 2.3
in a way different from Lemma 2.1.

\par \vskip 1.5pc
Let  $S_c$  be a smooth rational surface with a relatively minimal
elliptic fibration  $\psi : S_c \rightarrow {\bold P}^1.$
Suppose that  $\psi$  has two fibers  ${\overline F}_1, \,$ 
$m_2{\overline F}_2$ such that 
either one of the following two cases occurs.

\par \vskip 1.5pc
Case(a).  $m_2 = 1, \,$ $\psi$  is multiple-fiber free, and
each  ${\overline F}_i$  is either one of
Kodaira types II, \, III, \, IV \, and $I_{n_i}$ ($n_i \ge 1$).

\par
Case(b). $m_2 = 2, \,$ ${\overline F}_2$  is 
of Kodaira type  $I_{s_2}$ ($s_2 \ge 0$), the multiplicity-two fiber  
$2{\overline F}_2$  is the only multiple fiber of  $\psi,$
and  ${\overline F}_1$  is either one of
Kodaira types II, \, III, \, IV \, and $I_{n_1}.$

\par \vskip 1.5pc
In Case(a) (resp. Case(b)), let  $g : S \rightarrow S_c$  be the
blowing-up of intersection points in  ${\overline F}_i$  for  
$i = 1$  and  $2$ (resp. for  $i = 1$  only)  and 
their infinitely near points so that  $F_i := g^*{\overline F}_i$
fits one of Cases ($\alpha$), ($\beta$), ($\gamma$) and
($\delta$) in Lemma 1.5, according to the type of  ${\overline F}_i.$
We shall use the notation  $Supp F_i = \sum_{j=1}^{n_i} D_{i,j}
+ \sum_j H_{i,j}$  there.
In Case(b) we let  $F_2 := g^*{\overline F}_2,$
where  $F_2$ ($\cong {\overline F}_2$)  is of Kodaira type  $I_{s_2}.$

\par
In Case(a) and Case(b), the canonical divisor formula
implies respectively the following two relations:
$$(2.4.1a) \,\,\,\,
{\Cal O}(-K_{S_c})^{\otimes 2} \cong 
{\Cal O}({\overline F}_1 + {\overline F}_2),$$
$$(2.4.1b) \,\,\,\,
{\Cal O}(-K_{S_c})^{\otimes 2} \cong {\Cal O}({\overline F}_1).$$
Now (2.4.1a) and (2.4.1b) induce respectively 
the following two relations:
$$(2.4.2a) \,\,\,\, {\Cal O}(-K_S)^{\otimes 2} \cong 
{\Cal O}(\sum_{i=1}^2 \sum_{j=1}^{n_i} D_{i,j}),$$
$$(2.4.2b) \,\,\,\, {\Cal O}(-K_S)^{\otimes 2} \cong 
{\Cal O}(\sum_{j=1}^{n_1} D_{1,j}).$$

\par
In Case(a) (resp. Case(b)), let  $\pi_c : X_c = 
Spec \oplus_{i=0}^1 {\Cal O}(iK_{S_c}) \rightarrow S_c$  
be the double covering associated with the relation
(2.4.1a) (resp. (2.4.1b)), and  let  
$\pi : X = Spec \oplus_{i=0}^1 {\Cal O}(iK_S) \rightarrow S$  
be the double covering associated with the relation
(2.4.2a) (resp. (2.4.2b)).  Then both  $Gal(X_c/S_c)$  and  
$Gal(X/S)$  are isomorphic to the same cyclic group  $<\sigma>$  
of order  $2$  such that  $X_c/\sigma = S_c$  and  $X/\sigma = S.$

\par \vskip 1.5pc
{\bf Lemma 2.4.} {\it Let  $S_c$  be a smooth rational surface 
with a relatively minimal elliptic fibration
$\psi : S_c \rightarrow {\bold P}^1$  fitting the above Case$($a$)$ 
$($resp. Case$($b$))$.  Let  $\pi_c : X_c \rightarrow S_c = X_c/\sigma$
and  $\pi : X \rightarrow S = X/\sigma$  be as above.
Then the following four asssertions hold true.}

\par
(1) {\it For both  $i = 1$  and  $2$ $($resp. $i = 1$ only$)$,
one has  $\pi^*D_{i,j} = 2C_{i,j}$  
for a smooth rational curve  $C_{i,j}$ $(j = 1,2,\cdots,n_i)$,
and the fixed locus  $X^{\sigma}$  is a disjoint union of
$n := n_1+n_2$ $($resp. $n:= n_1)$ curves  $C_{i,j}$  contained in  
``the fiber''  $E_i$ {\rm (see (4) below)}.  One has  
the Picard number  $\rho(X) \ge \rho(S) = 10 + n,$  and
$2 \le n \le 10$ $($resp. $1 \le n \le 9)$.}

\par
(2) {\it $X$  is a smooth  K3 surface.  The involution  $\sigma$
on  $X$  satisfies  $\sigma^*\omega = -\omega$  for
a non-zero holomorphic 2-form  $\omega$  on  $X.$}

\par
(3) {\it There exists a  $\sigma$-equivariant birational 
morphism  $f : X \rightarrow X_c$  which induces the following
commutative diagram with  $f$  as the minimal resolution and 
$g$  a resolution}
$$X \overset{f} \to{\rightarrow} X_c$$
$$\pi \downarrow \hskip2pc \downarrow \pi_c$$
$$S \overset{g} \to{\rightarrow} S_c.$$

\par
(4) {\it For both  $i = 1$  and  $2$ $($resp. $i = 1$ only$)$, 
one has  $\pi^*F_i = 2E_i$  where   $E_i$  
fits one of Cases {\rm (A), (B), (C)} and {\rm (D)} in Lemma $1.5$
according to the type of  $F_i.$
In Case$($b$)$, one has  $\pi^*F_2 = E_2$  where  $E_2$  is
of Kodaira type  $I_{2s_2}.$}

\par \vskip 1.5pc
{\it Proof.} The first three assertions, except for the 
last part of (1), can be proved similarly as in Lemma 2.1, 
while the assertion(4) follows from the construction of  $\pi.$

\par
Note that  $\rho(S) = 10 - K_S^2 = 10 + n$  by (2.4.2).

\par
When  ${\overline F}_i$  is of Kodaira type II, III or IV
one has  $n_i = 1, 2$  or  $4$  respectively.  Now assume
that  ${\overline F}_i$  has Kodaira type  $I_{n_i}.$
Since  $\psi : S_c \rightarrow {\bold P}^1$  is relatively
minimal, one has  $K_{S_c}^2 = 0$  and the Picard number
$\rho(S_c) = 10.$  Now the inequalities for  $n$  follow from
the inequality  $\rho(S_c) \ge 2 + \sum_i 
(\#$(irreducible components in  ${\overline F}_i$) $-1$).
This proves Lemma 2.4.

\par \vskip 1.5pc
In the sense of Horikawa [H1, p.48], we make the following:

\par \vskip 1.5pc
{\bf Definition 2.5.}  In Lemma 2.4, the double covering
$\pi : X \rightarrow S$  is called {\it the canonical resolution}
of the covering  $\pi_c : X_c \rightarrow S_c.$

\par \vskip 1.5pc
{\bf Remark 2.6.} The pair  $(n_1,n_2)$  in Lemma 2.4 Case(a), satisfies:
$$n_i \ge 1, \,\, 2 \le n_1+n_2 \le 10.$$
On the other hand, [P, the list] (see also [Mi]) classified 
the set of singular fibers of a rational relatively minimal elliptic fibration
$\psi : S_c \rightarrow {\bold P}^1$  with no multiple fibers.
Going through Persson's list, we see that all pairs  $(n_1,n_2)$  
satisfying the above conditions, can be obtained in the way of Lemma 2.4.
In particular, $n := n_1+n_2$ in Lemma 2.4 Case(a), can
attain any value in the range  $2 \le n \le 10.$

\par \vskip 1.5pc
{\bf Example 2.7.}  Let  ${\overline D} \subseteq {\bold P}^1$  
be a rational curve of degree  $6$  with  $10$  ordinary nodes
as its only singular points.  Let  $h : S \rightarrow {\bold P}^1$
be the blowing-up of all these  $10$  nodes and let  
$D = h'{\overline D}$  be the proper transform.
Then one has 
${\Cal O}(-K_S)^{\otimes 2} \cong {\Cal O}(D).$

\par
Denote by  $q : S \rightarrow R$  the contraction of
$D$  into a cyclic quotient singular
point of Brieskorn type  $C_{4,1}.$  Then one has
the relation  ${\Cal O}(-K_R)^{\otimes 2} \cong {\Cal O}_R.$
a rational log Enriques surface of index 
$2$  with exactly one singular point.

\par \vskip 1.5pc
{\bf Example 2.8.} By [P, the list], there is a smooth rational
surface  $S_c$  with a multiple-fiber free relatively minimal
elliptic fibration  $\psi : S_c \rightarrow {\bold P}^1$
such that  $\psi$  has two singular fibers  ${\overline F}_i$
($i = 1,2$) of Kodaira type  $I_{n_i},$
where  $(n_1,n_2)$  can be taken as any of the following
three pairs:  $(1,9), \, (2,8), (5,5)$ (cf. Lemma 7 in \S 7).
In particular, we have  $n_1+n_2 = 10.$

\par
We use the notation in Lemma 2.4 Case(a):
$g : S \rightarrow S_c, \,\, F_i = g^*{\overline F}_i =
D_{i,1}+2H_{i,1}+D_{i,2}+2H_{i,2}+ 
\cdots + D_{i,n_i-1}+2H_{i,n_i-1}+D_{i,n_i}, \,\,$
and the canonical resolution  $\pi : X \rightarrow S = X/\sigma$
of the double covering  $\pi_c : X_c \rightarrow S_c.$

\par
Let  $M'$  be a  $(-1)$-curve on  $S_c.$  The relation (2.4.1a)
implies that  $M' . {\overline F}_i = 1.$  Let  $M = g^{-1}M',$
which is a  $(-1)$-curve on  $S.$  We may assume that
$M. D_{1,n_1} = M . D_{2,1} = 1,$  after relabelling.
Thus we obtain a linear chain of  $10$
$(-4)$-curves  $D_{i,j}$  and  $9$  $(-1)$-curves
$M$  and  $H_{i,j}$  as follows:
$$D_{1,1}-H_{1,1}-D_{1,2}-H_{1,2}- \cdots -D_{1,n_1}-
M-D_{2,1}-H_{2,1}-D_{2,2}-H_{2,2}- \cdots -D_{2,n_2}.$$

\par
Actually, in view of Lemmas 2.1(4), 2.3 and 2.4, starting 
with any two fibers  ${\overline F}_i$
of Kodaira type II (setting  $n_i := 1$), III ($n_i := 2$), 
IV ($n_i := 4$)  or  $I_{n_i}$  with  $n_1+n_2 = 10,$ 
we can obtain a linear chain of  $10$  $(-4)$-curves 
and  $9$  $(-1)$-curves as shown above,
after a suitable relabelling.

\par \vskip 1.5pc
{\bf (2.8.1).}  Let  $v_2 : S \rightarrow S_2$  be the 
contraction of this linear chain of length  $19$  into 
a cyclic quotient singular point of Brieskorn type  $C_{40, 19}$ 
(Lemma 1.4).  Then  $S_2$  is isomorphic to the {\it unique} 
rational log Enriques surface of Type  $A_{19}$ (see 
[OZ1, Theorem 2 and Example 2] and [Z2, Example 3.2] for the 
proof of the uniqueness and two different constructions, of
the same surface).  By Definition 2.2,
our  $\pi$  is also the canonical resolution of the
canonical covering  ${\overline \pi}_2 : X_2 \rightarrow S_2.$

\par \vskip 1.5pc
{\bf (2.8.2).}  Let  $v_3 : S \rightarrow S_3$  be the 
contraction of  $10$  disjoint  $(-4)$-curves  $D_{i,j}$  
into  $10$  cyclic-quotient singular points 
of Brieskorn type  $C_{4,1}.$  Then  $S_3$  is the 
{\it unique} rational log Enriques surface of index  $2$ 
with  $10$  Cartier-index two singular points (Corollary 5).  
We note that  $\sum_{i,j} D_{i,j} = \pi(X^{\sigma})$ (Lemma 2.4).

\par \vskip 1.5pc
{\bf (2.8.3).} Since  $X^{\sigma}$  is a disjoint union
of  $10$  smooth rational curves  $C_{i,j} := \pi^{-1}(D_{i,j}),$
the pair  $(X, \sigma)$  here is isomorphic to Shioda-Inose's 
unique pair in Theorem 3.  In particular,
$S = X/\sigma$  is independent of the choice of  
fibers  ${\overline F}_i$  so long as  $n_1+n_2 = 10.$

\par \vskip 1.5pc
{\bf Lemma 2.9.} {\it Let  $(X, \sigma)$  be as in Theorem $1$.
Assume further that  $X^{\sigma}$  is a disjoint union of
$n$ $(n \ge 1)$  smooth rational curves  $C_i.$}

\par
{\it Then the pair  $(X, \sigma)$  can be realized in the way 
of Lemma $2.4$, from a rational surface  $S_c$  satisfying
all hypotheses there, i.e., our quotient morphism  
$\pi : X \rightarrow S := X/\sigma$  is the canonical resolution 
of the double covering  $\pi_c : X_c \rightarrow S_c$  there.
Moreover, each fiber  ${\overline F}_i$  there can be so chosen
that it is not of Kodaira type $IV$.}

\par \vskip 1.5pc
{\it Proof.} By Lemma 2.3, the quotient morphism  
$\pi : X \rightarrow S := X/\sigma$   is the canonical resolution,
in the sense of Definition 2.2 or Lemma 2.1,
of the canonical covering  ${\overline \pi} : W \rightarrow R$
of a rational log Enriques surface  $R$  of index  $2$
with a cyclic-quotient singular point of Brieskorn type 
$C_{4n, 2n-1}$  as its only singular point.

\par
We use the notation in Lemma 2.1 :  $q_1 : S_1 \rightarrow R, \,
q_2 : S \rightarrow S_1, \, q = q_1 \circ q_2, \,$
and the commutative diagram there.  Note that  
$\Delta := q^{-1}(Sing R) = \sum_{i=1}^n D_i + \sum_{j=1}^{n-1} H_j$  
is a linear chain as shown in the proof of Lemma 2.1(4); 
$X_2, \, S_2$  there and  $W, \, R$  there coincide 
in the present case.  $q_2$  is just the smooth blowing-down 
of  $n-1$  $(-1)$-curves  $H_j$  in  $\Delta.$

\par
By [GZ3] or [Z6, Proposition 3.1], there exists a  
$(-1)$-curve  $M'$  on  $S_1$  such that  $M'$  meets a 
tip component of the linear chain  
$q_2(\Delta).$  Note that  $M' . q_2(\Delta) = 2$
because of the relation  $-2K_{S_1} \sim q_2(\Delta)$  
deduced from (2.1.2).
Let  $M := q_2'M'$  be the proper transform on  $S.$  
Then one of the following three cases occurs.

\par \vskip 1.5pc
Case($\alpha$)'.  $M^2 = -1, \,$ $M . \Delta = 2\,$  and  $M$  
has an order-two touch with  $D_1$  at a point  $d_1 \in D_1 - D_2.$ 
In this case we set  $n_1 := 1.$

\par
Case($\beta$)'. $M^2 = -2,\,$ $M . \Delta = 1,$  and  $M$  
meets  $H_1$  at a point  $d_1 \in H_1 - (D_1+D_2).$  
In this case we set  $n_1 := 2.$

\par
Case($\delta$)'.  $M^2 = -1, \,$ $M . \Delta = 2$  and  $M$  
meets two distinct points  $d_1, d_2$  of  $\sum_i D_i.$
To be precise, $d_1 \in D_1, \, d_2 \in D_{n_1}$ ($n_1 \ge 1$).

\par \vskip 1.5pc
In Case ($\alpha$)', ($\beta$)' and ($\delta$)', 
let  $F_1 := 2M+D_1, \,$  $F_1 := 4H_1 + 2M + D_1 + D_2,$
$F_1 := D_1 + 2H_1 + D_2 + 2H_2 + \cdots + D_{n_1-1} + 2H_{n_1-1}
+ D_{n_1} + 2M,$  respectively.
Then  $F_1$  fits respectively Case ($\alpha$), or ($\beta$), or ($\delta$)
of Lemma 1.5, so that  $F_1$  contains  $\sum_{i=1}^{n_1} D_i.$

\par \vskip 1.5pc
{\bf Claim(1).}  There is an elliptic fibration 
$\psi : S \rightarrow {\bold P}^1$  with  $F_1$  as
a non-multiple fiber.

\par \vskip 1.5pc
By the construction of  $\pi,$  one has  $\pi^*F_1 = 2E_1$  where
$E_1$  fits Case (A), (B) or (D) in Lemma 1.5 according to
the type of  $F_1.$  By the Riemann-Roch theorem, there exists
an elliptic fibration  $\varphi : X \rightarrow {\bold P}^1$ 
with  $E_1$  as a fiber.  Clearly,  $X^{\sigma} = \sum_{i=1}^{n_1} C_i,$
where  $\pi^*D_i = 2C_i,$  is contained in fibers 
of  $\varphi$  and  $E_1$  is  $\sigma$-stable.  So one can
apply Lemma 1.11.  Thus Claim(1) follows.

\par \vskip 1.5pc
{\bf Claim(2).}  Let  $E_1, \, E_2$  be the only  $\sigma$-stable fibers of
$\varphi$ (Lemma 1.11).  Then  $E_2$  does not fit Case(C) in Lemma 1.5.

\par \vskip 1.5pc
If  $n_1 = n$  then  $E_2$  contains no component of
$X^{\sigma}$  and hence  $E_2$  is of Kodaira type  $I_{s_2}.$
If  $n > n_1,$ all components
of  $E_2 \cap X^{\sigma}$  are contained in the linear
chain  $\pi^{-1}(D_{n_1+1}+H_{n_1+1}+\cdots+D_n)$
which is a subset of  $E_2.$  Then Claim(2) is clear
because there is no such a linear chain in  $E_2$  fitting Case(C).

\par \vskip 1.5pc
By Lemma 1.5 and Claim(2), $E_2$ fits Case (a) or (b) below.
Let  $g : S \rightarrow S_c$  be the smooth blowing-down
of curves in fibers  $\pi_*E_i$  so that  
$\psi : S \rightarrow {\bold P}^1$  induces a relatively
minimal elliptic fibration, also denoted by  
$\psi : S_c \rightarrow {\bold P}^1$ (Lemma 1.11 (3)).
Then  ${\overline F}_1 := g_*F_1$  has respectively
Kodaira type (II), (III) or  $I_{n_1}.$ 

\par \vskip 1.5pc
Case(a). $E_2$  fits one of Cases (A), (B) and (D) in Lemma 1.5.
Then  $F_2 := \pi_*E_2$  is a non-multiple fiber of  $\psi$
fitting respectively one of Cases ($\alpha$), ($\beta$) and 
($\delta$) in Lemma 1.5.  Hence  $\psi$  is multiple fiber 
free (Lemma 1.11 (3)).  Moreover, ${\overline F}_2 := g_*F_2$
has respectively one of  Kodaira types II, III and  $I_{n_2}$ 
($n_2 \ge 1$), where  $n = n_1 + n_2.$

\par
Case(b). $E_2$ is of Kodaira type  $I_{2s_2}$ ($s_2 \ge 0$).
To be precise, $E_2$  is either an elliptic curve or fits 
Case(E) in Lemma 1.5.  So  $\pi_*E_2 = 2F_2$  where  $F_2$ 
is of Kodaira type  $I_{s_2}.$  Hence  $2F_2$  or  
$2{\overline F}_2,$  where  ${\overline F}_2 := g_*F_2$ 
($\cong F_2$), is the only multiple fiber of  $\psi$ 
(Lemma 1.11 (3)).  One has  $n = n_1.$

\par \vskip 1.5pc
Thus  $S_c$  satisfies all hypotheses in Lemma 2.4.
By Lemma 1.2 (2), our quotient morphism  $\pi : X \rightarrow
S = X/\sigma$  is the canonical resolution of the
double covering  ${\overline \pi}_c : X_c \rightarrow S_c$
associated with the relation (2.4.1a) or (2.4.1b), respectively,
because our blowing-down  $g : S \rightarrow S_c$
coincides with the map  $g : S \rightarrow S_c$  
constructed preceding Lemma 2.4.  This proves Lemma 2.9.

\par \vskip 1.5pc
{\bf Definition and Proposition 2.10.} In Lemma 2.9, 
the double covering  $\pi_c : X_c \rightarrow S_c$  is called 
{\it a canonical mapping model} of  $\pi : X \rightarrow S.$
This  $\pi_c$  is not unique (see Example 2.8), for
there is no unique elliptic fibration  
$\varphi : X \rightarrow {\bold P}^1$  such that  
$X^{\sigma}$  is contained in fibers of  $\varphi.$
One may have also noticed that the fiber  
${\overline F}_i$  in Lemma 2.4 can be taken as
of Kodaira type IV, but we can avoid this type in Lemma 2.9.

\par \vskip 2pc
{\bf \S 3. Log del Pezzo surfaces of Cartier index}  ${\bold \le 2}$

\par \vskip 1pc
Let  $S_c$  be a log del Pezzo surface of Cartier-index
$\le 2.$  We shall use the notation in Lemma 1.4, where  $r \ge 0:$
$$g_1 : S_1 \rightarrow S_c (= T), \,\,
g_2 : S \rightarrow S_1, \,\, g = g_1 \circ g_2, \,\,
K_S = g^*K_{S_c} - 1/2 \sum_{i=1}^r \sum_{j=1}^{n_i} D_{i,j}.$$

\par \vskip 1.5pc
{\bf Lemma 3.1.} {\it Let  $S_c$  be a log del Pezzo surface
of Cartier-index  $\le 2.$  Then the following two assertions 
are true.}

\par
(1) {\it $|-2K_{S_c}|$  contains a smooth curve  
${\overline F}$  of genus  $\ge 2$
with  ${\overline F} \cap Sing S_c = \phi.$}

\par
(2) {\it {\rm dim} $|-2K_{S_c}| = 3K_{S_c}^2 = 3(n + K_S^2),$
where  $n := \sum_{i=1}^r n_i.$}

\par \vskip 1.5pc
{\it Proof.}  (1) When  $S_c$  has Cartier-index  $2$, (1)
is proved in [AN, Theorem 3].
 
\par
Now assume that  $S_c$  is a log del Pezzo surface of 
Cartier-index one, i.e., $S_c$  has at worst
Du Val singular points.  
Then  $-K_S = -g^*K_{S_c}$  is nef and big, and
has zero intersection with  $g^{-1}(Sing S_c).$
By [D, Theorem 1, p. 55], $|-2K_S|$  is base point free.
Actually, [D] assumed the condition that  $S$  is a blowing-up of
several points on  ${\bold P}^2.$  However, if this condition
is not true then  $S$  is the Hirzebruch surface of degree  $2$  and  
$g : S \rightarrow S_c$  is the contraction of the
$(-2)$-curve.  In this case,
$|-2K_{S_c}|$  is also base point free.

\par
So a general member  $F$  of  $|-2K_S|$  is smooth and also 
connected because  $(-2K_S)^2 > 0.$  Now  ${\overline F} := g_*F$  
satisfies the conditions in (1).  Indeed, 
$g({\overline F}) = g(F) \ge 2$  comes from the genus formula.
This proves (1).

\par
(2) Set  $F := g^*{\overline F} \sim g^*(-2K_{S_c}).$
Since  $S_c$  has only rational singularities and by 
the projection formula, one has  $R^ig_*{\Cal O}(F) \cong
{\Cal O}({\overline F}) \otimes R^ig_*{\Cal O}_S = 0$
for all  $i \ge 0.$  Hence  $H^i(S, {\Cal O}(F)) \cong$
$H^i(S_c, g_*{\Cal O}(F)) = H^i(S_c, {\Cal O}({\overline F}))$
for all  $i.$  Since ${\overline F} - K_{S_c} \sim -3K_{S_c}$
is  ${\bold Q}$-ample, [KMM, Theorem 1-2-5] implies that
$H^i(S_c, {\overline F}) = 0$  for all  $i \ge 1.$  

\par
Now the Riemann-Roch theorem implies that  
$h^0(S_c, -2K_{S_c}) = h^0(S, F) = 1 + 1/2g^*$
$(-2K_{S_c}) . (g^*(-2K_{S_c}) - K_S)$
$= 1 + 1/2g^*(-2K_{S_c}) . (g^*(-2K_{S_c}) - g^*K_{S_c})$
$= 1 + 3K_{S_c}^2.$  By the ramification formula preceding
Lemma 3.1, one has  $K_S^2 = K_{S_c}^2 - n.$
Now (2) follows.  This completes the proof of Lemma 3.1.

\par \vskip 1.5pc
Let  ${\overline F}$  be as in Lemma 3.1.
Then one has:
$$(3.2.1) \,\,\,\,\,
{\Cal O}(-K_{S_c})^{\otimes 2} \cong {\Cal O}({\overline F}).$$
Let  $F := g^*{\overline F},$  which is a smooth curve
away from  $g^{-1}(Sing S_c).$  Now (3.2.1) induces:
$$(3.2.2) \,\,\,\,\, -2K_S \sim 
g^*(-2K_{S_c}) + \sum_{i=1}^r \sum_{j=1}^{n_i} D_{i,j}
\sim F + \sum_{i,j} D_{i,j}, \,\,
{\Cal O}(-K_S)^{\otimes 2} \cong {\Cal O}(F + \sum_{i,j} D_{i,j}).$$
Here  $F + \sum_{i,j} D_{i,j}$  is a disjoint union of
$F$  and  $n$  $(-4)$-curves  $D_{i,j},$  where  $n = \sum_{i=1}^r n_i.$
Let  
$\pi_c : X_c := Spec \oplus_{i=0}^1 {\Cal O}(iK_{S_c}) \rightarrow S_c$
and  
$\pi : X := Spec \oplus_{i=0}^1 {\Cal O}(iK_S) \rightarrow S$
be the double coverings associated with the relations (3.2.1)
and (3.2.2), respectively.
Then both  $Gal(X_c/S_c)$  and  $Gal(X/S)$  are 
isomorphic to the same cyclic group  $<\sigma>$  of order  $2$  
such that  $X_c/\sigma = S_c$  and  $X/\sigma = S.$

\par \vskip 1.5pc
{\bf Lemma 3.2.} {\it Let  $S_c$  be a log del Pezzo surface
of Cartier-index  $\le 2.$  Let  
$\pi_c : X_c \rightarrow S_c = X_c/\sigma$  and  
$\pi : X \rightarrow S$  be as above.
Then the following three assertions hold true.}

\par
(1) {\it One has  $\pi^*F = 2E$  for a smooth curve  $E$  isomorphic
to  $F,$  and  $\pi^*D_{i,j} = 2C_{i,j}$  for a smooth rational
curve  $C_{i,j}.$  The fixed locus  $X^{\sigma}$  is 
a disjoint union of  $E$  and  $n$  curves  $C_{i,j},$  
where  $n := \sum_{i=1}^r n_i.$  One has  $n \le 9.$}

\par
(2) {\it $X$  is a smooth  K3 surface.  The involution  $\sigma$
on  $X$  satisfies  $\sigma^*\omega = -\omega$  for
a non-zero holomorphic $2$-form  $\omega$  on  $X.$}

\par
(3) {\it There exists a  $\sigma$-equivariant birational 
morphism  $f : X \rightarrow X_c$  which induces the following 
commutative diagram with  $f$  as the minimal resolution
and  $g$  a resolution}
$$X \overset{f} \to{\rightarrow} X_c$$
$$\pi \downarrow \hskip2pc \downarrow \pi_c$$
$$S \overset{g} \to{\rightarrow} S_c.$$

\par
(4) {\it Suppose in addition that  $n = 9.$  Then  $S_c$  is
isomorphic to the unique log del Pezzo surface of Picard 
number  $1$  and Cartier index  $2$ {\rm (see \S 7 for the
construction and [AN, Figure 1] for
the configuration of all  $19$  exceptional curves on  $S$)}.}

\par
{\it In particular, $Sing S_c$  is a single cyclic quotient 
singular point of Brieskorn type  $C_{36, 17},\,$
$g(E) = 2$  and the Picard number  $\rho(S) = 18,$ 
whence  $\rho(X) \ge 18.$  Moreover, 
${\Cal O}(X^{\sigma}) = {\Cal O}(E + \sum_{i,j} D_{i,j})$
is divisible by  $2$  in  $Pic S$  and hence  $S$  is the 
quotient of a smooth surface of general type, modulo an
involution.}

\par \vskip 1.5pc
{\it Proof.}  Except for the inequality  $n \le 9,$ assertions
(1), (2) and (3) can be proved similarly as in Lemma 2.1.
If  $S_c$  has Cartier index  $1$  then  $n = 0.$

\par
If  $S_c$  has Cartier index  $2,$  then  $1 \le n \le 9$
by [AN, Theorem 4].  Suppose that  $n = 9.$  Then (4)
follows from [AN, Theorems 4 $\sim 7$], because there is apparently
no proper DPN-subgraph  $\Gamma(S)$  of  $\Gamma({\widetilde S}),$
where  ${\widetilde S}$  is extremal with  $(n, g(E), \delta) = (9,2,0)$
in the notation there, due to the fact that  
$DuVal \Gamma({\widetilde S}) = \phi.$

\par \vskip 1.5pc
{\bf Definition 3.3.}  In Lemma 3.2, the double covering
$\pi : X \rightarrow S$  is called {\it the canonical resolution}
of the covering  $\pi_c : X_c \rightarrow S_c.$

\par \vskip 1.5pc
The following Lemma 3.4 is the converse to Lemma 3.2.

\par \vskip 1.5pc
{\bf Lemma 3.4.} {\it Let  $(X, \sigma)$  be as in Theorem 1.
Assume further that  $X^{\sigma}$  contains a smooth curve
$E$  of genus  $\ge 2.$  

\par
Then  $X^{\sigma}$  is a disjoint union
of  $E$  and  $n$ $(n \ge 0)$ smooth rational curves,
and the pair  $(X, \sigma)$  can be realized in the way of Lemma $3.2$,
from a log del Pezzo surface  $S_c$  of Cartier index  $\le 2,$
i.e., the quotient morphism  $\pi : X \rightarrow S := X/\sigma$
is the canonical resolution of the double covering
$\pi_c : X_c \rightarrow S_c$  there.}

\par \vskip 1.5pc
{\it Proof.}  Let  $\Sigma$  be the set of all
curves on  $X$  which have zero intersection with  $E.$
Since  $E^2 = 2g(E) -2 > 0,$  the Hodge index theorem
and the finiteness of  rank $Pic X$
imply that  $\Sigma$  consists of finitely many curves
and has negative definite intersection matrix.
In particular, every curve in  $\Sigma$  is 
a smooth rational curve by the genus formula.
Thus every connected component of  $\Sigma$
is disjoint from  $E$  and has Dynkin type  $A_m$ 
($m \ge 1$), $\, D_m$ ($m \ge 4$) or  $E_m$ ($m = 6,7,8$).

\par
Let  $\Gamma_i$ ($i = 1,2,\cdots,r; r \ge 0$)  be all connected 
components of  $\Sigma$  containing a curve of  $X^{\sigma}.$
Let  $\Delta_j$ ($j = 1,2, \cdots, s_1; s_1 \ge 0$)  
be the remaining connected components of  $\Sigma.$

\par
Since  $E$  is  $\sigma$-fixed, $\sigma$  induces a 
permutation on the set of connected components of  $\Sigma.$
So  $\Gamma_i$  is  $\sigma$-stable,
while  $\sigma(\Delta_j) = \Delta_{j'}$  for some  $j' \ne j$
(Lemma 1.6).  Hence we may assume that  $s_1 = 2s$
and  $\sigma(\Delta_i) = \Delta_{s+i}$ ($i = 1,2,\cdots,s$).

\par
On the other hand, by Lemma 1.6, $\Gamma_i$  is of type 
$A_{2n_i-1}$ ($n_i \ge 1$)  as follows:
$$C_{i,1}-G_{i,1}-C_{i,2}-G_{i,2}- \cdots 
-C_{i,n_i-1}-G_{i,n_i-1}-C_{i,n_i}.$$
Here  $C_{i,j}$  is  $\sigma$-fixed, while  $G_{i,j}$  is
$\sigma$-stable but not  $\sigma$-fixed.

\par
Note that  $X^{\sigma}$  is the disjoint union of
$E$  and  $n$  smooth rational curves  $C_{i,j},$ where  
$n = \sum_{i=1}^r n_i.$  So, by Lemma 1.2, the quotient 
morphism  $\pi : X \rightarrow S = X/\sigma$ 
coincides with the double covering associated with
the relation  (3.2.2)  where  $F := \pi_*E = \pi(E), \,$
$D_{i,j} := \pi_*C_{i,j} = \pi(C_{i,j}).$

\par
Since  $\pi^*F = 2E,$ the following is clear.

\par \vskip 1.5pc
{\bf Claim(1).} $\pi(\Sigma)$  consists of exactly all 
curves on  $S$  having zero intersection with  $F.$

\par \vskip 1.5pc
Clearly,  $\pi(\Delta_j) = \pi(\Delta_{s+j})$ ($j = 1,2,\cdots,s$)
is a connected component of  $\pi(\Sigma)$  disjoint from
$F$  and with the same weighted 
dual graph as  $\Delta_j.$

\par
On the other hand, $\pi(\Gamma_i)$ ($i = 1,2,\cdots,r$) 
is a connected component of  $\pi(\Sigma)$  disjoint from
$F$  and with the following dual graph :
$$D_{i,1}-H_{i,1}-D_{i,2}-H_{i,2}- 
\cdots -D_{i,n_i-1}-H_{i,n_i-1}-D_{i,n_i}.$$
Here  $D_{i,j}$  is a  $(-4)$-curve, while  
$H_{i,j}: = \pi(G_{i,j})$  is a  $(-1)$-curve.

\par
Let  $g : S \rightarrow S_c$  be the contraction of  $\pi(\Sigma)$
into points.  Then  $t_i := g\pi(\Gamma_i)$ ($i = 1,2,\cdots,r$) 
is a cyclic quotient singularity of Brieskorn type  
$C_{4n_i, 2n_i-1}$ (Lemma 1.4), 
while  $g\pi(\Delta_j)$  is a Du Val singular point.

\par
By Claim(1), ${\overline F} := g_*F$  is a smooth ample 
Cartier divisor isomorphic to  $F$ (and also to  $E$).
Our relation (3.2.2) induces the relation (3.2.1).
Thus, $S_c$  satisfies the hypothesis in Lemma 3.2.

\par
It is clear that  $\pi : X \rightarrow S$  is the
canonical resolution of the double covering  
$\pi_c : X_c \rightarrow S_c$  associated with the relation 
(3.2.1), because our contraction  $g : S \rightarrow S_c$  
here coincides with the map  $g : S \rightarrow S_c (=T)$  
constructed in Lemma 1.4.  This proves Lemma 3.4.

\par \vskip 1.5pc
{\bf Definition  and Proposition 3.5.}  In Lemma 3.4, 
the double covering  $\pi_c : X_c \rightarrow S_c$  is called 
{\it the canonical mapping model} of  $\pi : X \rightarrow S.$
This  $\pi_c$  is unique, because the map  $g : S \rightarrow S_c$  
is the contraction of all curves on  $S$  having zero intersection 
with  $\pi_*E$  and uniquely determined by  $\pi.$ 

\par \vskip 2pc
{\bf \S 4.  Multiple-fiber free rational elliptic fibrations}

\par \vskip 1pc
Let  $S$  be a smooth rational surface 
with a multiple-fiber free relatively minimal elliptic fibration
$\psi : S \rightarrow {\bold P}^1.$
Let  $F_1, \, F_2$  be two smooth fibers of  $\psi.$
By the canonical divisor formula, one has
$K_S \sim -F_1.$  Hence one obtains the following relation:
$$(4.1.1) \,\,\,\,
{\Cal O}(-K_S)^{\otimes 2} \cong {\Cal O}(F_1+F_2).$$

Let  $\pi : X := Spec \oplus_{i=0}^1 {\Cal O}(iK_S) \rightarrow S$
be the double covering associated with the relation (4.1.1).
Then  $Gal(X/S)$  is a cyclic group  $<\sigma>$  of
order  $2$  such that  $X/\sigma = S.$

\par
The following lemma can be proved similarly as in Lemma 2.1.

\par \vskip 1.5pc
{\bf Lemma 4.1.} {\it Let  $S$  be a smooth rational surface 
with a multiple-fiber free relatively minimal elliptic fibration
$\psi : S \rightarrow {\bold P}^1.$  
Let  $\pi : X \rightarrow S = X/\sigma$  be as above.
Then the following two asssertions hold true.}

\par
(1) {\it One has  $\pi^*F_i = 2E_i$  for an elliptic
curve  $E_i$  isomorphic to  $F_i.$  The fixed locus
$X^{\sigma}$  is a disjoint union of  $E_1$  and  $E_2.$}

\par
(2) {\it $X$  is a smooth  K$3$ surface.  The involution  $\sigma$
on  $X$  satisfies  $\sigma^*\omega = -\omega$  for
a non-zero holomorphic $2$-form  $\omega$  on  $X.$}

\par \vskip 1.5pc
The following Lemma 4.2 is the converse to Lemma 4.1.

\par \vskip 1.5pc
{\bf Lemma 4.2.} {\it Let  $(X, \sigma)$  be as in Theorem $1$.
Assume further that  $X^{\sigma}$  is a disjoint union of two 
elliptic curves  $E_1, \, E_2.$  Then the pair  $(X, \sigma)$  
can be realized in the way of Lemma $4.1$, from a rational surface
$S$  satisfying all hypotheses there.}

\par \vskip 1.5pc
{\it Proof.}  By Lemma 1.11, $S = X/\sigma$  is a smooth rational
surface and there exists a multiple-fiber free relatively
minimal elliptic fibration  $\psi : S \rightarrow {\bold P}^1$  
with  $F_i := \pi_*E_i = \pi(E_i)$  as smooth fibers.
Now Lemma 4.2 follows from Lemma 1.2 (2).

\par \vskip 2pc
{\bf \S 5. Rational elliptic fibrations with a multiple fiber}

\par \vskip 1pc
Let  $S$  be a smooth rational surface 
with a relatively minimal elliptic fibration  
$\psi : S \rightarrow {\bold P}^1$  such that  $2F_2$  
is the only multiple fiber of
$\psi,$  where  $F_2$  is of Kodaira type  $I_s$ ($s \ge 0$).

\par
By the canonical divisor formula, one has
$K_S \sim -F_1 + F_2$  for a smooth fiber  $F_1$
of  $\psi.$  Hence one obtains the following:
$$(5.1.1) \,\,\,\, 
{\Cal O}(-K_S)^{\otimes 2} \cong {\Cal O}(F_1).$$

Let  $\pi : X := Spec \oplus_{i=0}^1 {\Cal O}(iK_S) \rightarrow S$
be the double covering associated with the relation (5.1.1).
Then  $Gal(X/S)$  is a cyclic group  $<\sigma>$  of
order  $2$  such that  $X/\sigma = S.$

\par
The following Lemma 5.1 can be proved similarly as in Lemma 2.1.

\par \vskip 1.5pc
{\bf Lemma 5.1.} {\it Let  $S$  be a smooth rational surface 
with a relatively minimal elliptic fibration
$\psi : S \rightarrow {\bold P}^1$  such 
that  $2F_2$  is the only multiple fiber of
$\psi,$  where  $F_2$  is of Kodaira type  $I_s$ $(s \ge 0)$.
Let  $\pi : X \rightarrow S = X/\sigma$  be as above.
Then the following two asssertions hold true.}

\par
(1) {\it One has  $\pi^*F_1 = 2E_1$  for an elliptic
curve  $E_1$  isomorphic to  $F_1.$  The fixed locus
$X^{\sigma}$  is equal to  $E_1.$}

\par
(2) {\it $X$  is a smooth  K$3$ surface.  The involution  $\sigma$
on  $X$  satisfies  $\sigma^*\omega = -\omega$  for
a non-zero holomorphic 2-form  $\omega$  on  $X.$}

\par \vskip 1.5pc
The following Lemma 5.2 is the converse to Lemma 5.1.

\par \vskip 1.5pc
{\bf Lemma 5.2.} {\it Let  $(X, \sigma)$  be as in Theorem $1$.
Assume further that  $X^{\sigma}$  is a single 
elliptic curve  $E_1.$  Then the pair  $(X, \sigma)$  can
be realized in the way of Lemma $5.1$, from a rational surface
$S$  satisfying all hypotheses there.}
  
\par \vskip 1.5pc
{\it Proof.} By Lemma 1.11, there exists a fiber  $E_2$
of  $\varphi := \Phi_{|E_1|} : X \rightarrow {\bold P}^1$
such that  $E_1, \, E_2$  are only  $\sigma$-stable
fibers of  $\varphi.$  Applying Lemma 1.5 to  $E_2,$  we see 
that  $E_2$  is of Kodaira type $I_{2s}$ ($s \ge 0$). 

\par
By Lemmas 1.11 and 1.5, $S = X/\sigma$  is a smooth rational surface with
a relatively minimal elliptic fibration  
$\psi : S \rightarrow {\bold P}^1$  such that  $F_1 := \pi_*E_1$  
is a smooth fiber of  $\psi$  and  $\pi_*E_2 = 2F_2$  is the only 
multiple fiber of  $\psi,$  where  $F_2$  is of Kodaira
type  $I_s.$  Now Lemma 5.2 follows from Lemma 1.2 (2).
This completes the proof of Lemma 5.2.

\par \vskip 2pc
{\bf \S 6. Multiple-fiber free rational elliptic fibrations with 
a fiber of Kodaira type II, III, IV  or}  ${\bold I_n}$

\par \vskip 1pc
Let  $S_c$  be a smooth rational surface 
with a multiple-fiber free relatively minimal 
elliptic fibration  $\psi : S_c \rightarrow {\bold P}^1.$
Suppose that  $\psi$  has a singular fiber  
${\overline F}_{\infty}$  of either one of the 
Kodaira types II, III, IV and  $I_{n_{\infty}}$ 
($n_{\infty} \ge 1$).

\par
By the canonical divisor formula, One has
$K_{S_c} \sim -{\overline F}_1$  for a smooth fiber  
${\overline F}_1$  of  $\psi.$
Hence one obtains the following relation:

$$(6.1.1) \,\,\,\,\, {\Cal O}(-K_{S_c})^{\otimes 2} \cong 
{\Cal O}({\overline F}_1 + {\overline F}_{\infty}).$$

\par
Let  $g : S \rightarrow S_c$  be the composite of blowing-ups 
of intersections of  ${\overline F}_{\infty}$  and their infinitely
near points so that  $F_{\infty} := g^*{\overline F}_{\infty}$
fits respectively one of Cases ($\alpha$), ($\beta$), ($\gamma$) 
and ($\delta$) in Lemma 1.5.  We shall use the
notation  $Supp F_{\infty} =$ 
$\sum_{j=1}^{n_{\infty}} D_{{\infty},j} + \sum_j H_{{\infty},j}$
there.  Here  $\sum_j D_{{\infty},j}$  is a disjoint union
of  $n_{\infty}$  $(-4)$-curves  $D_{{\infty},j}.$  

\par
The relation (6.1.1) induces the following relation, where 
$F_1 := g^*{\overline F}_1$
$$(6.1.2) \,\,\,\,\,  
{\Cal O}(-K_S)^{\otimes 2} \cong 
{\Cal O}(F_1 + \sum_{j=1}^{n_{\infty}} D_{{\infty},j}).$$

Let  
$\pi_c : X_c:= Spec \oplus_{i=0}^1 {\Cal O}(iK_{S_c}) \rightarrow S_c$
and  $\pi: X := Spec \oplus_{i=0}^1 {\Cal O}(iK_S) \rightarrow S$
be the double coverings associated with the relations (6.1.1)
and (6.1.2), respectively.
Then both  $Gal(X_c/S_c)$  and  $Gal(X/S)$  are isomorphic 
to a cyclic group  $<\sigma>$  of order  $2$  
such that  $X_c/\sigma = S_c$  and  $X/\sigma = S.$

\par
The following Lemma 6.1 can be proved similarly as in Lemma 2.4.
In fact, the second part of the assertion(4)
follows from the observation:
the Picard number  $\rho(S) = 10 - K_S^2 = 10 + n_{\infty}$
by (6.1.2).

\par \vskip 1.5pc
{\bf Lemma 6.1.} {\it Let  $S_c$  be a smooth rational surface 
with a multiple-fiber free relatively minimal elliptic fibration
$\psi : S_c \rightarrow {\bold P}^1.$  Assume further that   
$\psi$  has a singular fiber  ${\overline F}_{\infty}$  of
Kodaira type II, III, IV or  $I_{n_{\infty}}$ $(n_{\infty} \ge 1)$.
Let  $\pi_c : X_c \rightarrow S_c = X_c/\sigma$
and  $\pi : X \rightarrow S = X/\sigma$  be as above.  
Then the following four asssertions hold true.}

\par
(1) {\it One has  $\pi^*F_1 = 2E_1$  for an elliptic
curve  $E_1$  isomorphic to  $F_1,$  and  
$\pi^*D_{{\infty},j} = 2C_{{\infty},j}$  for a smooth 
rational curve  $C_{{\infty},j}.$  The fixed locus  
$X^{\sigma}$  is a disjoint union of  $E_1,$  and  
$n_{\infty}$  curves  $C_{{\infty},j}$  all contained in  
``the fiber''  $E_{\infty}$ {\rm (see (4) below)}.  
One has  $n_{\infty} \le 9.$}

\par
(2) {\it $X$  is a smooth  K$3$ surface.  The involution  $\sigma$
on  $X$  satisfies  $\sigma^*\omega = -\omega$  for
a non-zero holomorphic $2$-form  $\omega$  on  $X.$}

\par
(3) {\it There exists a  $\sigma$-equivariant birational 
morphism  $f : X \rightarrow X_c$  which induces the following
commutative diagram with  $f$  as the minimal resolution and 
$g$  a resolution}
$$X \overset{f} \to{\rightarrow} X_c$$
$$\pi \downarrow \hskip2pc \downarrow \pi_c$$
$$S \overset{g} \to{\rightarrow} S_c.$$

\par
(4) {\it One has  $\pi^*F_{\infty} = 2E_{\infty},$  
where   $E_{\infty}$  fits one of Cases {\rm (A), (B), (C)} and {\rm (D)} 
in Lemma $1.5$ according to the type of  $F_{\infty}.$
One has also  $\rho(X) \ge \rho(S) = 10 + n_{\infty}.$}

\par \vskip 1.5pc
{\bf Definition 6.2.} In Lemma 6.1, the double covering  
$\pi : X \rightarrow S$  is called {\it the canonical resolution} 
of  $\pi_c : X_c \rightarrow S_c.$

\par \vskip 1.5pc
The following Lemma 6.3 is the converse to Lemma 6.1.

\par \vskip 1.5pc
{\bf Lemma 6.3.} {\it Let  $(X, \sigma)$  be as in Theorem $1$.
Assume further that  $X^{\sigma}$  is a disjoint union of
an elliptic curve  $E_1$  and  $n$ $(n \ge 1)$  smooth rational
curves.}  

\par
{\it Then the pair  $(X, \sigma)$  can be realized in the way 
of Lemma $6.1$, from a rational surface  $S_c$  satisfying 
all hypotheses there, i.e., the quotient morphism  
$\pi : X \rightarrow S := X/\sigma$
is the canonical resolution of the double covering
$\pi_c : X_c \rightarrow S_c$  there.}

\par \vskip 1.5pc
{\it Proof.} By Lemma 1.11, $\varphi := \Phi_{|E_1|} :$ 
$X \rightarrow {\bold P}^1$  has exactly two  $\sigma$-stable fibers
$E_1, \, E_{\infty}.$  Applying Lemma 1.5, $E_{\infty}$  fits one of
Cases (A), (B), (C) and (D) there with  $n_{\infty} = n,$  
and  $F_{\infty} := \pi_*E_{\infty}$
fits respectively one of Cases ($\alpha$), ($\beta$), ($\gamma$)
and ($\delta$) there.  We use the notation 
$Supp E_{\infty} = \sum_{j=1}^n C_{{\infty},j} + \sum_j G_{{\infty},j}$  
and  
$Supp F_{\infty} = \sum_{j=1}^n D_{{\infty},j} + \sum_j H_{{\infty},j}$  
there, where  $\pi^* D_{{\infty},j} = 2C_{{\infty},j}.$
Now Lemmas 1.11 and 1.5 imply that  $S$  is a smooth 
rational surface, which has a multiple-fiber free elliptic 
fibration  $\psi : S \rightarrow {\bold P}^1$  with  
$F_1 := \pi_*(E_1) = \pi(E_1)$  as its smooth fiber.

\par
Let  $g : S \rightarrow S_c$  be the smooth blowing-down of
curves in  $F_{\infty}$  so that  
${\overline F}_{\infty} := g_*F_{\infty}$  
is of Kodaira type II, III, IV or 
$I_n,$ according to the type of  $E_{\infty}.$  Now  
$\psi : S \rightarrow {\bold P}^1$  induces a 
multiple-fiber free relatively minimal elliptic fibration 
(Lemmas 1.11), also denoted by  
$\psi : S_c \rightarrow {\bold P}^1.$  This  $S_c$  
satisfies all hypotheses of Lemma 6.1.  

\par
By Lemma 1.2 (2), our quotient morphism 
$\pi : X \rightarrow S = X/\sigma$  coincides with
the double covering associated with the relation (6.1.2).
Thus  $\pi$  is the canonical resolution of the double covering
$\pi_c : X_c \rightarrow S_s$  associated with the relation 
(6.1.1), because our contraction 
$g : S \rightarrow S_c$  here coincides with the map  
$g : S \rightarrow S_c$  constructed at the beginning 
of \S 6.  This proves Lemma 6.3.

\par \vskip 1.5pc
{\bf Definition and Proposition 6.4.} 
In Lemma 6.3, the double covering  $\pi_c : X_c \rightarrow S_c$  
is called  {\it the canonical mapping model} of  $\pi : X \rightarrow S.$
This  $\pi_c$  is unique, because  $\varphi : X \rightarrow {\bold P}^1$
is the only elliptic fibration such that  $X^{\sigma}$  is 
contained in fibers of  $\varphi,$  and there is a unique 
smooth blowing-down  $g : S \rightarrow S_c$  of curves in  
$F_{\infty}$  such that  $g_*F_{\infty}$  is a minimal fiber.

\par \vskip 2pc
{\bf \S 7. Proofs of Theorems 1, 3 and 4 and Corollary 5}

\par \vskip 1pc
Theorem 1 is a consequence of Lemmas 1.2, 1.11, 2.9, 3.4, 
4.2, 5.2 and 6.3.  Theorem 4 follows from Lemmas 2.1 and 2.9.

\par \vskip 1.5pc
Next we prove Corollary 5 using Theorem 3'.  
The first part of (1) is proved in Lemma 2.1(4).
By Remark 2.6 and Example 2.7, for any  $1 \le n \le 10,$  
there is a K3 surface  $X$  with an involution  $\sigma$  
such that  $X^{\sigma}$  is a disjoint union of  $n$
smooth rational curves  $C_i.$  

\par
Let  $\pi : X \rightarrow S := X/\sigma$  be the
quotient morphism.  As in Lemma 2.1 (4), let  
$v_3 : S \rightarrow S_3$  be the contraction of
$n$  $(-4)$-curves  $D_i := \pi_*C_i.$  Then  $S_3$
is a rational log Enriques surface of index  $2$  with  
$n$  cyclic-quotient singular points  $v_3(D_i)$  of 
Brieskorn type  $C_{4,1}$  as its only singular points
(Lemma 1.4).  This proves Corollary 5(1).

\par
Now let  $R$  be a rational log Enriques surface of
index  $2$  with  $10$  singular points  $r_i$  of Cartier 
index  $2.$  Let  $\pi : X \rightarrow S := X/\sigma$  
be the canonical resolution of the canonical covering  
${\overline \pi} : W \rightarrow R$ (Definition 2.2).
By Lemma 2.1, one has  $r = 10$  and  $n = \sum_{i=1}^{10} n_i
\le 10.$  So, $n_i = 1$  for all  $i.$  
Hence  $Sing R$  consists of  $10$  cyclic-quotient 
singular points  $r_i$  of Brieskorn type  $C_{4,1}$
and several Du Val singular points.

\par
In the notation of Lemma 2.1, $D_i := q^{-1}(r_i)$
is a  $(-4)$-curve on  $S$  and  $X^{\sigma}$  
is the disjoint union of  $10$  smooth
rational curves  $C_i$  where  $\pi^*D_i = 2C_i.$
By Theorem 3,  $(X, \sigma)$  is
isomorphic to Shioda-Inose's unique pair.
Now we have only to show that  $R$  has no Du Val singular 
points, for then  $q : S = X/\sigma \rightarrow R$  is 
just the contraction of  $\pi(X^{\sigma})$  and Corollary 5(2)
follows from the uniqueness of the pair  $(X, \sigma).$

\par
Suppose to the contrary that  $R$  has a Du Val singular 
point  $r_0.$  Then  $q^{-1}(r_0)$  consists of  $(-2)$-curves
disjoint from  $\sum_{i=1}^{10} D_i.$
Hence  $\pi^{-1}q^{-1}(r_0)$  consists of smooth rational curves
disjoint from  $\sum_{i=1}^{10} C_i = X^{\sigma}.$
By the following Claim(1), each  $(-2)$-curve on
$X$  is  $\sigma$-stable and we reach a contradiction
to Lemma 1.1 (3).  

\par \vskip 1.5pc 
{\bf Claim(1).}  $\sigma^*|Pic X = id.$

\par \vskip 1.5pc
This is proved in [OZ1, Lemma 3.3].  Actually, in the notation 
of Lemma 2.1 (4), $\sigma$  stabilizes  $19$  curves  
$\sum_{i=1}^{10} C_i + \sum_{j=1}^9 G_j$
as well as the pull back of the generator of  $Pic X_2.$

\par
This proves Claim(1) and also Corollary 5.

\par \vskip 1.5pc
Now we prove Theorem 3' in \S 1 which is stronger than Theorem 3.  
The inequality  $0 \le m \le 10$
follows from Theorem 1.  One sees that  $m$  can attain any
value in this range by Theorem 1(1) and the proof of Corollary 5(1).
Assume that  $m = 10.$  Then, by Theorem 1, 
the quotient morphism  $\pi : X \rightarrow S := X/\sigma$  
is the canonical resolution of its canonical mapping model  
$\pi_c : X_c \rightarrow S_c$  given in either 
Lemma 2.4 with  $n = 10,$ or Lemma 3.2 
with  $n = 9,$ or Lemma 6.1 with  $n_{\infty} = 9$  and 
${\overline F}_{\infty}$  of Kodaira type  $I_9.$ 
Now Theorem 3'(1) follows. 

\par
We now prove Theorem 3'(2).  Then  $\pi_c$  is given 
in Lemma 2.4 with  $n = 10.$  By Lemma 2.3, our
$\pi : X \rightarrow S = X/\sigma$  
is the canonical resolution of the canonical double 
covering  ${\overline \pi} : W \rightarrow R$  of a 
rational log Enriques surface  $R$  of Type  $A_{19}.$  
This  $R$  determines uniquely the pair
$(X, \sigma)$  (see Lemma 2.1).  Now Theorem 3'(2)
follows from [OZ1, Theorem 2] saying that
there is, up to isomorphisms, only one 
rational log Enriques surface of Type  $A_{19}.$ 

\par
We need the following Lemma 7 to prove Theorem 3' (3)(4).

\par \vskip 1.5pc
{\bf Lemma 7.} {\it Let  $S_{ci}$ $(i = 1,2)$  be a 
smooth rational surface with a multiple-fiber free 
relatively minimal elliptic fibration  
$\psi_i : S_{ci} \rightarrow {\bold P}^1$
which has exactly one singular fiber  ${\overline F}_{ci}$  
of Kodaira type  $I_9$  and three singular fibers of 
Kodaira type  $I_1$  as its only singular fibers.
Let  $M_i$  be a  $(-1)$-curve on  $S_{ci}.$}

\par
{\it Then there is an isomorphism  
$\tau : S_{c1} \rightarrow S_{c2}$  such that  
$\tau^*{\overline F}_{c2}
= {\overline F}_{c1}$  and  $\tau^*M_2 = M_1.$}

\par
{\it Finally, there exists a smooth rational curve
$H_0$  on  $S_{c1}$  such  that  $H_0^2 = 0$  and  $H_0$
is a  $2$-section of  $\psi_1$  passing through
an intersection of  ${\overline F}_{c1}.$}

\par \vskip 1.5pc
{\it Proof.} By the hypothesis on  $\psi_i$  and
the canonical divisor formula, one has  
$K_{S_{ci}} + {\overline F}_{ci} \sim 0.$
Hence  $M_i$  is a cross-section.  Thus we can
write  ${\overline F}_{ci} = \sum_{j=1}^9 D_{i,j}$
so that  $M_i . D_{i,1} = D_{i,j} . D_{i,j+1} = 1,$
where  $D_{i,9+j} := D_{i,j}.$
Let  $h_i : S_{ci} \rightarrow {\bold P}^2$  be the
smooth blowing down of  $M_i+\sum_{j=1}^8 D_{i,j}$  into
the node  $d_i$  of the nodal cubic  $h_i(D_{i,9}).$

\par 
We may assume that both  $h_i(D_{i,9})$  are equal to
${\overline D} : Y^2Z - X^2(X+Z) = 0$  in  ${\bold P}^2.$  
Note that the projective transformation  
$\tau_1 : {\bold P}^2 \rightarrow {\bold P}^2,$
where  $\tau_1(X) = -X, \, \tau(Y) = Y, \, \tau(Z) = -Z,$
stabilizes  ${\overline D}$  and switchs two local irreducible 
components of  ${\overline D}$  at its node  $d := [0:0:1].$
Note also that  $\Lambda_i := h_{i*}|{\overline F}_{ci}|$  is a
$1$-dimensional linear system satisfying the following 
hypothesis($*$):

\par \vskip 1.5pc
(*) \hskip 2pc $\Lambda_i$  contains the nodal cubic  
${\overline D} : Y^2Z - X^2(X+Z) = 0$  as a member.
A general member of  $\Lambda_i$  is an elliptic curve
and touches, with order  $8,$  the local irreducible
component of  ${\overline D},$ tangent to  
$Y+\varepsilon_i X = 0,$  at the node  $d.$  
Here  $\varepsilon_i = \pm 1.$

\par \vskip 1.5pc
{\bf Claim(2).} For a given  $\varepsilon_i,$  the
$\Lambda_i$  above is the only  $1$-dimensional
linear system satisfying the hypothesis($*$) above.

\par \vskip 1.5pc
Indeed, suppose that  $G_1, G_2$  are two distinct 
elliptic curves, each of which touches, with order 8,
the local irreducible component of  ${\overline D}$,
tangent to  $Y + \varepsilon X = 0$, at the node  $d$.
Let  $h : W \rightarrow {\bold P}^2$  be the blowing-up
of  $d$  and 7 of its infinitely near points,
such that  ${\overline D}$  and the proper transform 
$G_i'$  on  $W$  of  $G_i$  have no intersection.  It is easy to see
that the pull-back on  $W$  of  ${\overline D}$
is a simple loop of  8  $(-2)$-curves and one  $(-1)$-curve
which together generate  $(Pic W) \otimes {\bold Q}$,
whence one can check that  $G_1', G_2'$  are numerically 
(and hence also linearly) equivalent.
Therefore, $G_1 \sim G_2$.  This proves Claim(1).

\par \vskip 1.5pc
There is a projective transformation  
$\tau : {\bold P}^2 \rightarrow {\bold P}^2$  such that  
$\tau^*\Lambda_2 = \Lambda_1.$  Indeed, let  $\tau = id$ 
if  $\varepsilon_1 = \varepsilon_2,$  and let  
$\tau = \tau_1$  above otherwise.  This  $\tau$  
induces an isomorphism between  $S_{c1}$  and  $S_{c2}$
required by Lemma 7, because the linear system
$|{\overline F}_i|$  is just the unique minimal resolution
of base points (i.e., $d_i$  and its infinitely near points) 
in  $\Lambda_i$  and hence uniquely determined by  $\Lambda_i.$

\par
For the last paragraph of Lemma 7, we take the projective
line  ${\overline H}_0$  through  $d_1$
and tangent to the local irreducible component,
other than the one in the hypothesis($*$), of  $h_1(D_{1,9})$  
at its node  $d_1.$  Then the proper transform
$H_0 := h_1'{\overline H}_0$  is through  $D_{1,8} \cap D_{1,9}$
and the one required by Lemma 7.  This completes the proof of Lemma 7.

\par \vskip 1.5pc
As a consequence to Lemma 7, one obtains:

\par \vskip 1.5pc
{\bf Corollary 8.} {\it There is, up to isomorphisms, only one
log del Pezzo surface  $T$  with a
type  $A_8$  Du Val singular point as its only singular point.}

\par \vskip 1.5pc
{\it Proof.}  Let  $T_i$  be two surfaces both satisfying all 
hypotheses of Corollary 8.  Let  $v_i : T_i' \rightarrow T_i$
be the minimal resolution.  Then  $\Delta_i := v_i^{-1}(Sing T_i)$  
has Dynkin type  $A_8.$  Note that  $-K_{T_i'} = -v_i^*K_{T_i}$
is nef and big.  Hence  $9 \ge 10 - (K_{T_i'})^2 = 
\rho(T_i') = 8 + \rho(T_i) \ge 9.$  So  $(K_{T_i'})^2 = 1$
and the Picard number  $\rho(T_i) = 1.$

\par
The Riemann-Roch theorem and the vanishing theorem 
[KMM, Theorem 1-2-3] imply that  dim $|-K_{T_i'}| = 1.$  
Thus, $Bs|-K_{T_i'}|$  consists of a single point  $t_i$  
and a general member  $F_i'$  of  $|-K_{T_i'}|$  
is an elliptic curve (see also [D, Theorem 1, p. 39]).

\par
Let  $w_i : S_{ci} \rightarrow T_i'$  be the blowing-up
of  $t_i.$  Then there exists a multiple-fiber free relatively
minimal elliptic fibration  $\psi_i : S_{ci} \rightarrow {\bold P}^1$
with  $w_i'F_i'$  as a general fiber
and the  $(-1)$-curve  $M_i := w_i^{-1}(t_i)$  as a cross-section.

\par
It is easy to see that the singular fiber  ${\overline F}_{ci}$  
of  $\psi_i$  containing  $w_i'\Delta_i$  is of Kodaira type  
$I_9.$  By [P, the list], our  $\psi$  satisfies all
hypotheses in Lemma 7.  So, there exists an isomorphism
$\tau : S_{c1} \rightarrow S_{c2}$  such that  $\tau^*M_2 = M_1$
and  $\tau^*{\overline F}_{c2} = {\overline F}_{c1}.$  
This  $\tau$  induces an isomorphism between  $T_1$  and  $T_2.$
This proves Corollary 8.

\par \vskip 1.5pc
We now continue the proof of Theorem 3' (3)(4).

\par
Let  $(S_c, \psi)$  be the unique pair, modulo isomorphisms,
in Lemma 7.  Set  ${\overline F}_s := \psi^{-1}(s).$  We may 
assume that  ${\overline F}_{\infty}, \, {\overline F}_0, \,
{\overline F}_1, \, {\overline F}_{s_0}$  are respectively
of Kodaira types  $I_9, I_1, I_1, I_1.$
Write  ${\overline F}_{\infty} = \sum_{i=1}^9 {\overline D}_i$  
so that  ${\overline D}_i . {\overline D}_{i+1} = 1,$
where  ${\overline D}_{9+i} := {\overline D}_i.$

\par
As in Lemma 6.1, let  $g_{ell} : S_{ell} \rightarrow S_c$  
be the smooth blowing-up of the  $9$  intersections in  
${\overline F}_{\infty}$  so that  
$F_{\infty} := g_{ell}^*{\overline F}_{\infty} =$
$D_1 + 2H_1 + D_2 + 2H_2 + \cdots + D_9 + 2H_9$  
with  $D_i := g_{ell}'{\overline D}_i,$  as given in 
Lemma 1.5 Case($\delta$).  $\psi$  on  $S_c$  induces
an elliptic fibration, also
denoted by  $\psi : S_{ell} \rightarrow {\bold P}^1$
with  $F_s := g_{ell}^*{\overline F}_s$  as a fiber.
Now a relation identical to (6.1.1) (but with different labelling) 
induces the following relation identical to (6.1.2)
except the labelling :
$$(3.2) \,\,\,\, {\Cal O}(-K_{S_{ell}})^{\otimes 2} \cong 
{\Cal O}(F_s + \sum_{i=1}^9 D_i).$$  

\par
As in Lemma 6.1, let  
$\pi : X_s \rightarrow S_{ell} = X_s/\sigma_s$  
be the double covering associated with the relation (3.2).
Then  $X_s/\sigma_s = S_{ell}$  is independent of 
the choice of  $s$  and Theorem 3'(3i) is true.

\par
When  $F_s$  is smooth, the pair  $(X_s, \sigma_s)$  
fits Lemma 6.1 with  $n_{\infty} = 9$  and is of 
Type(Ell).  Now Theorem 3'(iv) follows from 
Theorem 1 and Lemmas 6.1, 6.3 and 7.

\par
When  $s = \infty,$ the right hand side of the relation
(3.2) is divisible by  $2$  in  $Pic S_{ell}.$  
So Theorem 3'(3ii) is true.

\par
Suppose that  $s = 0, 1,$  or  $s_0.$  By the construction 
in Lemma 2.4 or Example 2.8, there is a smooth blowing-up 
$X/\sigma \rightarrow X_s/\sigma_s$  of the node 
of  $F_s,$ where  $(X, \sigma)$  is Shioda-Inose's unique 
pair, constructed also in Example 2.8 with  $(n_1,n_2=n_{\infty})$
$= (1,9).$  Thus Theorem 3'(3iii) is true.  This 
proves Theorem 3'(3).  

\par
By Lemma 7, there is a smooth rational
curve  ${\overline H}_0$  on  $S_c$  such that  
${\overline H}_0^2 = 0$  and that  ${\overline H}_0$  
is a  $2$-section of  $\psi$  through
the intersection  ${\overline D}_9 \cap {\overline D}_1$
after rotating the indices.  
Set  $H_0 := g_{ell}^{-1}({\overline H}_0),$  which is a
$(-1)$-curve.  For each of  $s = 0,1,s_0$, 
$H_0$  intersects the fiber  $F_s$ 
at its two distinct smooth points, for otherwise
the inverse on   $X$  of  $H_0$  is a 
union of a  $(-2)$-curve and its  $\sigma$-conjugate
and we reach a contradiction to the fact that
$\sigma^*|Pic X = id$ (see Claim(1) in the proof of Corollary 5).
Now since  $H_0$  is a 2-section of  $\psi$,
$H_0$  has a contact of order 2 with a smooth fiber
$F_{s_1}$  say, and intersects each fiber  $F_s$ ($s \ne \infty, s_1$)
at its two smooth points.

\par
Let  $S_{ell} \rightarrow S_{gn2}$  be the smooth blowing-down
of  $H_0.$  Then (3.2) induces the following
relation, where  $F_s', \, D_i'$  and also  $H_i'$
(for later use) are images of  $F_s, \, D_i, \, H_i$
$$(3.3) \,\,\,\, {\Cal O}(-K_{S_{gn2}})^{\otimes 2} \cong
{\Cal O}(F_s' + \sum_{i=1}^9 D_i').$$
Note that  $F_{\infty}'$  is a non-reduced simple loop.

\par
Let  $g_{gn2} : S_{gn2} \rightarrow {\overline S}_{gn2}$
be the contraction of the linear chain  
$D_1'+H_1'+D_2'+H_2' + \cdots +
D_9'$  into a cyclic quotient singularity of
Brieskorn type  $C_{36, 17}.$
Then one obtains:
$$(3.4) \,\,\,\, {\Cal O}(-K_{{\overline S}_{gn2}})^{\otimes 2} 
\cong {\Cal O}(g_{gn2*}F_s').$$

\par
Note that  $g_{gn2*}F_{\infty}'$  is twice of a nodal rational 
curve with the only singular point of  ${\overline S}_{gn2}$
as its node, and  $g_{gn2*}F_s'$ 
($s \in {\bold P}^1 - \{\infty\}$)  
is a curve of arithmetic genus  $2$  and
away from the singular point of  ${\overline S}_{gn2};$
moreover, for each of  $s = 0,1,s_0$  the curve  $g_{gn2*}F_s$
is rational with two simple nodes, $g_{gn2*}F_{s_1}$
is elliptic with a simple cusp, and for each  
$s \ne \infty,0,1,s_0,s_1$  the curve  $g_{gn2*}F_s$
is elliptic with a simple node.

\par
Since the Picard number  $\rho({\overline S}_{gn2}) =$
$\rho(S_{gn2}) -  17 = \rho(S_{ell}) - 18 = \rho(S_c) + 9 - 18 = 1,$
our  ${\overline S}_{gn2}$  is the unique log del Pezzo 
surface of Cartier index  $2$  and Picard number  $1$ 
(see Lemma 3.2(4)).

\par
By Lemmas 3.1 and 3.2(4), dim $|-2K_{{\overline S}_{gn2}}| = 3.$  
So  ${\bold P}(|-2K_{{\overline S}_{gn2}}|) =$ 
${\bold P}(|g_{gn2}^*(-2K_{{\overline S}_{gn2}})$
$|) = {\bold P}^3.$  
Let  $B$  denote the  $1$-dimensional linear subsystem
of  $|-2K_{{\overline S}_{gn2}}|$  consisting of
$g_{gn2*}F_s',$  the direct image of  
$F_s = g_{ell}^*{\overline F}_s$ ($s \in {\bold P}^1$).

\par
By the constructions, our  
$g_{gn2} : S_{gn2} \rightarrow {\overline S}_{gn2}$  
coincides with the  $g : S \rightarrow S_c$  in Lemma 1.4
or Lemma 3.2 with  $n = 9.$  As in Lemma 3.2, for any member  
$F_t' \in |g_{gn2}^*(-2K_{S_c})|$ ($t \in {\bold P}^3$), 
we let  $\pi : Y_t \rightarrow S_{gn2} = Y_t/\sigma_t$ 
be the double covering associated with the relation (3.3)
where  $F_s'$  is replaced by  $F_t'.$
Then  $Y_t/\sigma_t = S_{gn2}$  is independent of 
the choice of  $t$  and Theorem 3'(4i) is true.

\par
When  $F_t'$ ($t \in {\bold P}^3$)  is smooth, 
$\pi : Y_t \rightarrow S_{gn2} = Y_t/\sigma_t$  coincides 
with  $\pi : X \rightarrow S = X/\sigma$  given in
Lemma 3.2 with  $n = 9,$  and  $(Y_t, \sigma_t)$  is of 
Type(Gn2).  Conversely, by Theorem 1 and
Lemmas 3.2 and Lemma 3.4, every pair of Type(Gn2) is 
isomorphic to  $(Y_t, \sigma_t)$  for some  $t \in {\bold P}^3.$  
So Theorem 3'(4ii) is true.

\par
Theorem 3'(4iii) follows from the construction of
$S_{ell} \rightarrow S_{gn2}$  and the definition
of  $B$  above.  This completes the proof of Theorem 3'.

\par \vskip 2pc
{\bf References}

\par \vskip 1pc
[AN] V. A. Alexeev and V. V. Nikulin, Classification of del Pezzo
surfaces with log-terminal singularities of index  $\le 2,$
and involutions on K3 surfaces, {\it Soviet Math. Dokl.}
{\bf 39} (1989), 507 - 511.

\par
[Bl] M. Blache, The structure of l.c. surfaces of
Kodaira dimension zero, I, {\it J. Alg. Geom.}
{\bf 4} (1995), 137 - 179.

\par
[Br] E. Brieskorn, Rationale Singularit\"aten komplexer Fl\"achen,
{\it Invent. math.} {\bf 4} (1968), 336 - 358.

\par
[D] M. Demazure, Surfaces de del Pezzo - I, II, III, IV, V,
in : {\it Lecture Notes in Mathematics} {\bf 777} (1980), 22 - 69.

\par
[GZ1,2] R. V. Gurjar and D. -Q. Zhang, $\pi_1$  of smooth points
of a log del Pezzo surface is finite : I, II,
{\it J. Math. Sci. Univ. Tokyo,} {\bf 1} (1994), 137 - 180;
{\bf 2} (1995), 165--196.

\par 
[GZ3] R. V. Gurjar and D. -Q. Zhang, On the fundamental groups
of some open rational surfaces, {\it Math. Ann.}
{\bf 306} (1996), 15--30.

\par
[H1] E. Horikawa, On deformations of quintic surfaces,
{\it Inv. math.} {\bf 31}(1975), 43 - 85.

\par
[H2] E. Horikawa, Algebraic surfaces of general type with
small  ${\bold c}_1^2, III,$ {\it Invent. math.} 
{\bf 47} (1978), 209 - 248.

\par
[Ka1] Y. Kawamata, The cone of curves of algebraic varieties,
{\it Ann. of Math.} {\bf 119} (1984), 603 - 633.

\par
[KMM] Y. Kawamata, K. Matsuda and K. Matsuki, 
Introduction to the minimal model
problem, in : {\it Advanced Studies in Pure Mathematics,}
{\bf 10} (1987), pp. 283 - 360.

\par
[Ko] J. Koll\'ar, Flips and abundance for algebraic 
threefolds, {\it Ast\'erisque} {\bf 211} (1992).

\par
[Mi] R. Miranda, Persson's list of singular fibers
for a rational elliptic surfaces, {\it Math. Z.}
{\bf 205} (1990), 191 - 211.

\par
[Mo] D. R. Morrison, {\it On K3 surfaces with large Picard number,}
Invent. Math. {\bf 75} (1984), 105-121.

\par
[MZ1,2] M. Miyanishi and D. -Q. Zhang, Gorenstein log del 
Pezzo surfaces of rank one, I, II, 
{\it J. of Alg.} {\bf 118} (1988), 63 - 84; {\bf 156} (1993), 
183 - 193.

\par
[N1] V. V. Nikulin, Discrete reflections groups in
Lobachevsky spaces and algebraic surfaces,
In : {\it Proc. Internat. Congr. Math. (Berkeley, Calif. 1986),} Vol. 1,
Amer. Math. Soc. Providence, R.I. 1987, pp. 654 - 671.

\par
[N2] V. V. Nikulin, Factor groups of groups of automorphisms of
hyperbolic forms with respect to subgroups generated by
2-reflections, {\it J. Soviet Math.} {\bf 22} (1983), no. 4.

\par
[N3,4,5] V. V. Nikulin, Del Pezzo surfaces with log-terminal
singularities. I, II, III, {\it Math. USSR Sbornik,} 
{\bf 66} (1990), 231 - 248; {\it Math. USSR Izvestiya,} {\bf 33} (1989),
355 - 372;  {\it Math. USSR Izvestiya,} {\bf 35} (1990), 657 - 675.

\par
[OZ1,2] K. Oguiso and D. -Q. Zhang, On the most extremal 
log Enriques surfaces, I, II;
{\it Amer. J. Math.} {\bf 118} (1996), 1277-1297;
submitted 1996.

\par
[P] U. Persson, Configurations of Kodaira fibers on 
rational elliptic surfaces, {\it Math. Z.} {\bf 205} (1990),
1 - 47.

\par
[R] M. Reid, Campedelli versus Godeaux, In : {\it Problems in the 
theory of surfaces and their classification, Trento, October
1988, F. Catanese et al. ed.} Academic Press 1991, pp. 309 - 365.

\par
[S] F. Sakai, Anticanonical models of rational surfaces,
{\it Math. Ann.} {\bf 269} (1984), 389 - 410.

\par
[SI] T. Shioda and H. Inose, On singular K3 surfaces, in :
{\it Complex analysis and algebraic geometry,} Iwanami Shoten and
Cambridge University Press (1977), 119 - 136.

\par
[TY] H. Tokunaga and H. Yoshihara, Degree of irrationality
of Abelian surfaces, {\it J. of Alg.} {\bf 174} (1995), 1111 - 1121.

\par
[V] E. B. Vinberg, The two most algebraic K3 surfaces,
{\it Math.Ann.}{\bf 265}(1983),1 - 21.

\par
[Z1] D. -Q. Zhang, Logarithmic del Pezzo surfaces with
rational double and triple singular points, {\it Tohoku Math. J.}
{\bf 41} (1989), 399 - 452.

\par
[Z2,3] D. -Q. Zhang, Logarithmic Enriques surfaces, I, II, 
{\it J. Math. Kyoto Univ.} {\bf 31} (1991), 419 - 466; 
{\bf 33} (1993), 357 - 397.

\par
[Z4] D. -Q. Zhang, Algebraic surfaces with nef and big
anti-canonical divisor, {\it Math. Proc. Camb. Phil. Soc.}
{\bf 117} (1995), 161 - 163.

\par 
[Z5] D. -Q. Zhang, Algebraic surfaces with log canonical
singularities and the fundamental groups of their smooth parts,
{\it Transactions of A.M.S.} {\bf 348} (1996), 4175--4184.

\par
[Z6] D. -Q. Zhang, Normal algebraic surfaces with
trivial bicanonical divisor, {\it J. of Alg.} {\bf 186} (1996),
970--989.

\par \vskip 1.5pc
Department of Mathematics
\par 
National University of Singapore
\par
Lower Kent Ridge Road
\par
SINGAPORE 119260

\par
e-mail : matzdq$\@$math.nus.edu.sg

\enddocument

\enddocument